\documentclass[11pt]{article}
\usepackage{amsmath,amssymb,amsthm,bbm,natbib,color,hyperref}
\usepackage[margin=1in]{geometry}
\usepackage{graphicx}
\usepackage{float}
\usepackage{xcolor}

\usepackage{tikz}
\usetikzlibrary{arrows.meta,calc, patterns}

\newtheorem{theorem}{Theorem}
\newtheorem{lemma}{Lemma}
\newtheorem{proposition}{Proposition}

\newtheorem{definition}{Definition}

\newtheorem{assumption}{Assumption}

\newcommand{\R}{\mathbb{R}}

\newcommand{\EE}[1]{\mathbb{E}\left[{#1}\right]}

\newcommand{\PP}[1]{\mathbb{P}\left\{{#1}\right\}}

\newcommand{\eqd}{\stackrel{\textnormal{d}}{=}}
\newcommand{\One}[1]{{\mathbbm{1}}\left\{{#1}\right\}}
\newcommand{\one}[1]{{\mathbbm{1}}_{{#1}}}
\newcommand{\iidsim}{\stackrel{\textnormal{iid}}{\sim}}

\newcommand\independent{\protect\mathpalette{\protect\independenT}{\perp}}
\def\independenT#1#2{\mathrel{\rlap{$#1#2$}\mkern2mu{#1#2}}}

\newcommand{\cX}{\mathcal{X}}
\newcommand{\cY}{\mathcal{Y}}
\newcommand{\cT}{\mathcal{T}}
\newcommand{\bX}{\mathbf{X}}
\newcommand{\bY}{\mathbf{Y}}
\newcommand{\bx}{\mathbf{x}}
\newcommand{\by}{\mathbf{y}}
\newcommand{\dtv}{\mathsf{d}_{\mathsf{TV}}}

\makeatletter
\def\and{%
  \end{tabular}%
  \hskip .5em \@plus.17fil\relax
  \begin{tabular}[t]{c}}
\makeatother

\title{False positive control in time series coincidence detection}
\author{Ruiting Liang\thanks{NSF--Simons AI Institute for the Sky (SkAI), 172 E. Chestnut St., Chicago, IL 60611, USA}\,\,\thanks{Committee on Computational and Applied Mathematics, University of Chicago} \and Samuel Dyson\footnotemark[1]\,\,\thanks{Department of Physics, University of Chicago} \and Rina Foygel Barber\footnotemark[1]\,\,\thanks{Department of Statistics, University of Chicago} \and Daniel E. Holz\footnotemark[1]\,\,\thanks{Department of Physics, Department of Astronomy \& Astrophysics, Enrico Fermi Institute, Kavli Institute for Cosmological Physics, University of Chicago}}
\begin{document}
\maketitle

\begin{abstract}
    We study the problem of \emph{coincidence detection} in time series data, where we aim to determine whether the appearance of simultaneous or near-simultaneous events in two time series is indicative of some shared underlying signal or synchronicity, or might simply be due to random chance. This problem arises across many applications, such as astrophysics (e.g., detecting astrophysical events such as gravitational waves, with two or more detectors) and neuroscience (e.g., detecting synchronous firing patterns between two or more neurons). In this work, we consider methods based on \emph{time-shifting}, where the timeline of one data stream is randomly shifted relative to another, to mimic the types of coincidences that could occur by random chance. Our theoretical results establish rigorous finite-sample guarantees controlling the probability of false positives, under weak assumptions that allow for dependence within the time series data, providing reassurance that time-shifting methods are a reliable tool for inference in this setting. Empirical results with simulated and real data validate the strong performance of time-shifting methods in dependent-data settings.
\end{abstract}

\section{Introduction}

Coincidence detection over time arises in a wide range of domains, including astrophysics, neuroscience, and econometrics. Broadly speaking, the objective is to identify instances in which two or more processes exhibit temporally aligned or near-simultaneous activity beyond what would be expected by chance. These coincident patterns may indicate latent dependencies, shared external influences, or coordinated behaviors.

To make the problem more concrete, consider the setting of having two sources of data---say from two detectors---where for each one, we observe a data stream measured at times $t=1,\cdots,T$. We write $X_1,\dots,X_T$ to denote the measurements from the first stream, and $Y_1,\dots,Y_T$ for the second stream. For example, in the setting of measuring the activity of two neurons, $X_t$ and $Y_t$ may be binary variables indicating whether each neuron fired at time $t$; if the two neurons show more synchronous behavior than would be expected by random chance, this might be due to both neurons responding to the same input stimuli and/or connectivity between the two neurons. As another example, the detection of astrophysical events (such as flashes of electromagnetic radiation, or a gravitational wave signal) may be confirmed whenever transient signals are identified simultaneously in multiple detectors. Other applications where similar problems arise include healthcare (where we may have simultaneous measurements from multiple physiological sensors for a single patient), climate and air quality (for instance, identifying correlations between fluctuations in air quality and potential causal factors), and finance (where we may want to identify variables whose fluctuations are synchronized or are linked via a causal relationship).

In this paper, we are interested in problems of this type: namely, determining whether an apparent association between two data streams, $\bX=(X_1,\dots,X_T)$ and $\bY=(Y_1,\dots,Y_T)$, is stronger than what we might expect to see by random chance. This problem is challenging due to the complexity of working with time series data, which may exhibit features such as autocorrelation and other forms of dependence. At a high level, there are two different types of inference questions we might ask in these types of setting:
\begin{itemize}
    \item \emph{Detecting events}, where we ask whether a simultaneous burst of high signal from both data streams at time $t$ indicates that a real external event has occurred specifically at time $t$---e.g., if a gravitational wave is detected at time $t$, is this a true detection or a false positive? For this first type of inference question, a single detector would suffice, if we were confident that we had an accurate model for its noise properties (and could therefore determine whether a detected gravitational wave was, potentially, a likely false positive). But if we do not want to rely on potentially incorrect assumptions, we can instead rely on the fact that each detector can corroborate the findings of the other.
    \item \emph{Detecting synchronicity}, where we wish to test whether the two data streams appear to be correlated across the entire duration of time $t=1,\dots,T$---e.g., identifying connectivity between a pair of neurons. For this second type of inference problem, the actual question of interest relates to the relationship between the multiple data streams; we cannot pose a question about connectivity between neurons without measurements from multiple neurons.
\end{itemize}

\begin{figure}
    \centering
    \includegraphics[width=0.48\textwidth]{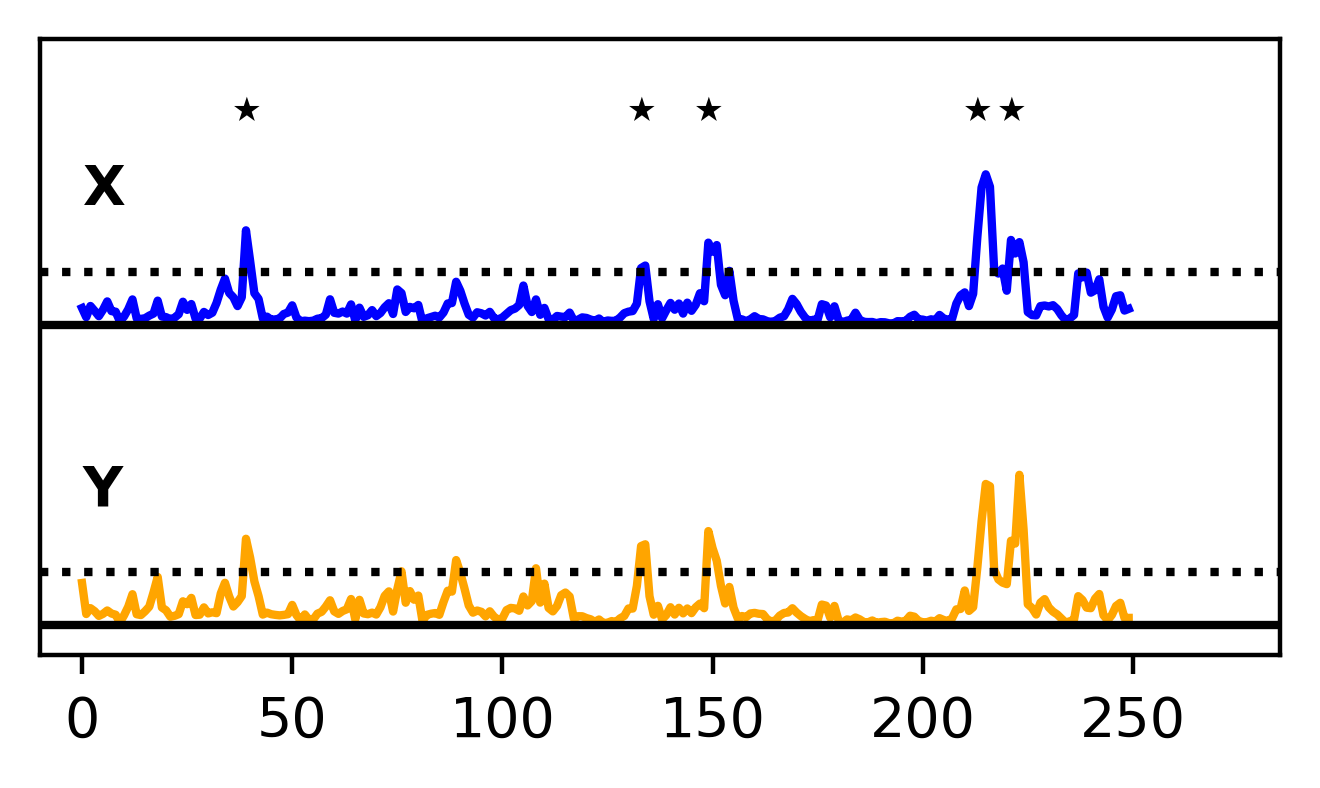}
    \quad 
    \includegraphics[width=0.48\textwidth]{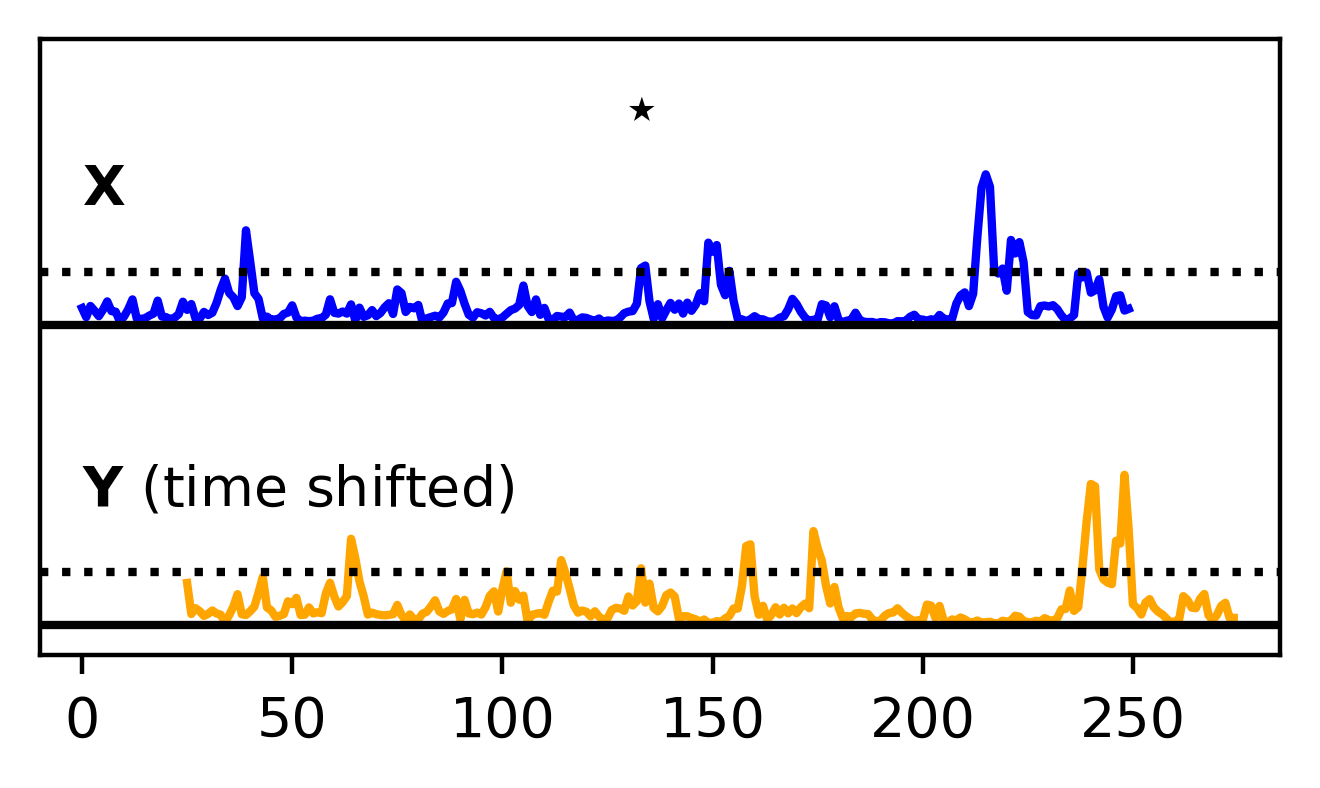}
    \caption{An illustration of the intuition behind the time-shifting approach. The left panel shows a dataset with time series $\bX$ and $\bY$, which exhibit some correlation. The star symbols in the figure mark times at which both $\bX$ and $\bY$ exceed a prespecified threshold (indicated by the dotted lines). On the right, we see a time-shifted version of the same dataset, where the $\bY$ time series has been time-shifted relative to $\bX$. We see that there are fewer stars (that is, there is less apparent association) within the time-shifted version of the data, suggesting that synchronicity may be present in the original data.}
    \label{fig:intro_timeshift}
\end{figure}

In many recent applications, questions of both of these types are often addressed using data-driven approaches, in order to avoid relying heavily on models that may not be perfectly accurate. A common strategy has emerged across many different fields, which relies (implicitly or explicitly) on the idea of \emph{time-shifting}, meaning that we analyze shifted versions of the two data streams. For example, if we are measuring synchronicity between the firing patterns of two neurons, we might compare the number of simultaneous spikes (the times $t$ where $X_t=1$ and $Y_t=1$), against a time-shifted version of the same statistic (the times $t$ where $X_t=1$ and $Y_{t-s}=1$, for some shift $s$). Effectively, when time-shifting, we are searching for synchronicity between the data streams $X_{1+s},X_{2+s},\dots$ and $Y_1,Y_2,\dots,\dots$, which allows us to determine what types of seemingly synchronous behavior might occur simply due to random chance, while preserving the complex properties of the data within each stream. See Figure~\ref{fig:intro_timeshift} for an illustration of the intuition behind this approach.

\subsection{Outline}
The remainder of this paper is organized as follows. In Section~\ref{sec:examples}, we explore several examples of applications of the coincidence detection problem, and describe some data-driven approaches that have been used for determining significance in each case. In Section~\ref{sec:method_events}, we will formally define the time-shifting method for detecting events, and will present theoretical guarantees on its performance, showing that the false detection rate of this type of procedure can be controlled with only mild assumptions. Section~\ref{sec:method_synchronicity} will do the same for the problem of detecting synchronicity.
We present empirical results on both simulated and real data in Section~\ref{sec:experiments}, and conclude with a discussion of related work, extensions and open questions for our findings, in Section~\ref{sec:discussion}. Proofs of all theoretical results are deferred to the Appendix.

\section{Motivating examples}\label{sec:examples}
In this section, we give a brief overview of several scientific contexts where the problem of coincidence detection has been addressed via a data-driven approach.

\paragraph{Synchronicity in neuron spike train data.}
In neuroscience, it is often important to determine whether multiple neurons are firing synchronously, which may indicate that they are interconnected and/or that they respond to the same stimuli. Suppose that two neurons are measured, with $X_i$ indicating whether the first neuron fired at time $i$, and similarly $Y_i$ indicating whether the second neuron fired at time $i$, for times $i=1,\dots,T$. The time series $X_1,\dots,X_T$ and $Y_1,\dots,Y_T$ are often referred to as the \emph{spike trains} of the two neurons. (Note that time has been discretized, i.e., each index $i$ corresponds to a short window of time, and $X_i$ or $Y_i$ indicates whether the neuron fired within that time window; this is the case for many applications.)

We may observe that there are multiple times $i$ when both neurons fired simultaneously---but is this higher than we would expect to see by random chance? A common approach towards answering this question is to examine the \emph{cross-correlation},
\begin{equation}\label{eq:Psi_crosscorr_neuro} \frac{1}{N_X N_Y}\sum_{i=1}^T \sum_{j=1}^T \One{X_i=1\textnormal{ and }Y_j=1\textnormal{ and }|i-j|\leq \Delta} ,\end{equation}
which measures the number of times that the two neurons fire at approximately the same time, across the entire time duration of measurement (with normalization by the total firing counts of each neuron, $N_X = \sum_{i=1}^T\one{X_i=1}$ and $N_Y = \sum_{i=1}^T \one{Y_i=1}$). Synchronicity measures of this type are common in the literature, e.g., \citet{perkel1967neuronal, Gerstein1972mutual, aertsen1989dynamics, EventsynchronizationQQ, hahnloser2002ultra}. 
Various data-driven approaches have been used to assess statistical significance of the resulting measure---for instance, \cite{pasquale2008self} describe a ``shuffling method'' that permutes time bins within each neuron's data stream to determine whether large values of the cross-correlation (or related measures) might occur by random chance. 

In other applications, the question of interest may be that of detecting events (e.g., was there a stimulus at a particular time that caused higher firing rates for both neurons?) rather than an overall pattern of synchronicity. In this case, the evidence at a given time $t$ can be computed as
\begin{equation}\label{eqn:Psi_event_neuro}\sum_{i=t}^{t+\Delta} \sum_{j=t}^{t+\Delta}X_i Y_j,\end{equation}
which counts only the number of simultaneous spikes from both neurons within a short time window after the stimulus occurrence at time $t$ (where $X_i,Y_i\in\{0,1\}$). Determining whether this measure of evidence is significant can again be done in a data-driven way, which is known as \emph{unitary event analysis} \citep{riehle1997spike}. In this setting, we might again have $X_i,Y_i\in\{0,1\}$ as before, or alternatively, if neurons are recorded over multiple trials (and the timing of stimuli occurring in the environment is the same across each trial), we might then use a statistic such as
$\sum_{i=t}^{t+\Delta} \sum_{j=t}^{t+\Delta} \sum_{k=1}^m X_{ik}  Y_{jk}$
to measure evidence of a response at time $t$, where $X_{ik},Y_{ik} \in\{0,1\}$ (with $X_{ik}\in\{0,1\}$ indicating whether the first neuron fired at time $i$ in trial $k$, and same for $Y_{ik}$).

\paragraph{Comparing flood events and epidemic outbreaks.}
While the previous examples consisted of two time series containing the same type of data (i.e., spike train data from two neurons), in this next example the two data streams are different in nature. \cite{donges2016event} study the question of whether flood events lead to epidemic outbreaks. Let $X_i$ be a binary variable indicating the presence of a flood event at time $i$, while $Y_i$ is binary and indicates the presence of an epidemic event. Since the hypothesis is that a flood event is a causal factor leading to an epidemic, the question is no longer symmetric: we are searching for flood events that occur shortly before (but not shortly after) epidemic events, and therefore we calculate
\[ \sum_{i=1}^T \One{\textnormal{$X_i=1$, and $Y_j=1$ for some $j$ with $i + a \leq j \leq i + b$}},\]
where $a<b$ are positive values indicating the window of time after a flood within which we might expect an epidemic to occur. 
After computing this test statistic,
in order to determine statistical significance,  \cite{donges2016event} propose an approximation to its null distribution (i.e., the distribution of the test statistic, in the absence of any true association between floods and epidemics), based on the observed number of floods and of epidemics. Concretely, the approximation of the null distribution is constructed by assuming that, in the absence of any association between epidemics and floods, each epidemic has probability $\frac{b-a}{T}$ of falling within the allowed time window $[i+a,i+b]$ occurring after a flood at time $i$.

\paragraph{Gravitational waves and gamma ray bursts.}
Multimessenger astronomy is the detection of transient events, such as a stellar explosion, using multiple types of detection signals (e.g., gamma-ray and optical observations). One such multimessenger detection occurred in 2017, with the near coincidence of space-based observations of a short gamma-ray burst (SGRB) 1.74 seconds after the ground-based detection of two merging neutron stars observed in gravitational waves (GWs) by the Laser Interferometer Gravitational-Wave Observatory (LIGO).

To place the coincidence detection question within the framework of our notation, let $X_i$ and $Y_i$ again be binary, with $X_i$ (respectively, $Y_i$) indicating the presence of a detection of a GW (respectively, a SGRB) at time $i$. Note that the question examines only a single GW detection---that is, $X_i=1$ for only a single time $i$ during the time window of observed data. The observed coincidence is quantified by the time gap between the GW detection, and the subsequent SGRB event that is closest in time:
\[ \min\{|i - t|: Y_i = 1\}\textnormal{ where }X_t=1,\]
i.e., where $t$ is the time of the GW detection. 
To evaluate the significance of the nearly simultaneous detection observed in 2017, \citet{Gravitational2017} uses the empirical rate of SGRB detections over multiple years of data, to evaluate the probability that these two events might have occurred nearly simultaneously by random chance---that is, this data-driven approach uses the SGRB data stream to estimate the frequency of events (i.e., how likely it is that $Y_i=1$, at each time $i$) in order to approximate the null distribution of the time gap test statistic. 

\paragraph{Multiple detectors for gravitational waves.}
Time stream data from LIGO's two detectors are used to make simultaneous, all-sky observations of gravitational waves (GWs), which are oscillations in spacetime produced by energetic events such as the merger of neutron stars and/or black holes \citep{AdvancedLIGO}.
Candidate events are identified through a template-based search in which each detector's data stream is compared against a pre-calculated signal waveform \citep{LIGO_Guide}. 
Let $\hat\rho_1(X_i)$ and $\hat\rho_2(Y_i)$ denote the signal-to-noise ratio (SNR) for each data stream, with higher values indicating greater alignment with the template. A simultaneous detection at time $t$ is then determined by the joint measure of evidence given by
\[ \max_{i : |i-t|\leq \Delta}\hat\rho(X_t,Y_i)\textnormal{ \ for \ }\hat\rho(X_t,Y_i) = \sqrt{\big(\hat\rho_1(X_t)\big)^2 + \big(\hat\rho_2(Y_i)\big)^2},\]
where $\Delta$ is the maximum possible time gap between the detections (determined by the distance between the two LIGO detectors). 
To determine whether a candidate event's observed measure of evidence might have simply occurred due to random chance, \cite{Usman_PyCBC} explicitly use a time-shifting approach, shifting the two data streams relative to each other and then re-calculating the values of $\hat\rho$---that is, for some shift $s$, we would consider
\[\hat\rho(X_t,Y_{i+s}),\]
for times $|i-t|\leq \Delta$.

\section{Coincidence detection via time-shifting: detecting events}\label{sec:method_events}
We now present a time-shifting method for the coincidence detection problem. In this section we focus on the first core question: the problem of detecting events (rather than detecting an overall pattern of synchronicity). However, as we will see later on, the two questions are closely related in terms of how time-shifting can be used for statistical inference.

\subsection{The method}

Suppose that, given our data streams $\bX = (X_1,\dots,X_T)$ and $\bY = (Y_1,\dots,Y_T)$ from two detectors, we would like to determine whether there is sufficient evidence from the detectors to detect an event that occurred at (approximately) time $t$. For instance, if these are measurements from two LIGO detectors, either one of the two detectors may exhibit a false positive (e.g., due to seismic vibrations from nearby vehicles), but if both detectors simultaneously identify a possible event then the two measurements corroborate each other and it is less likely to be a false positive. 

To formalize this setup, suppose that an event at time $t$ might plausibly affect the first detector during the time range $i\in\{t,\dots, t + \Delta\}$, and the second detector during the time range $j\in\{t,\dots,t + \Delta\}$, for some fixed and prespecified $\Delta\geq 0$. Define also some function $\psi$ that we will use as a measure of evidence---that is, the evidence for an event detection at time $t$ is captured by the quantity
\[\Psi_t := \psi\big( (X_t,\dots,X_{t+\Delta}),(Y_t,\dots,Y_{t+\Delta})\big)\in\R.\]
As a concrete example, we can consider the function defined by~\eqref{eqn:Psi_event_neuro}, in the context of unitary event analysis in neuroscience, where large values indicate a higher frequency of (nearly) simultaneous spikes in the two neurons. More generally, we will always take the convention that larger (positive) values indicate more evidence of a true detection.\footnote{In some settings, we may want to use an ``offset'' in time, for the two data streams---for instance, in the example of testing whether flood events trigger outbreaks of an epidemic, described in Section~\ref{sec:examples}. We could then use a test statistic of the form $\Psi_t = \psi( (X_{t+a_0},\dots,X_{t+a_1}), (Y_{t+b_0},\dots,Y_{t+b_1}))$, for some values $a_0,a_1,b_0,b_1$ that might differ for the two data streams. Our results can be extended in a straightforward way to this more general setting.}
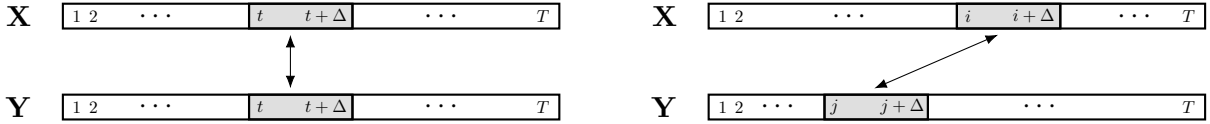
\begin{figure}
    \centering
        \begin{minipage}[t]{0.48\textwidth}
        \centering
        \begin{tikzpicture}[>=Latex, thick, scale=0.55]
        \def\xL{0}       
        \def\xR{12}      
        \def\barH{0.6}   
        \def\yGap{2.2}   
        \def\xLeftArrow{1.0}   
        \def\xRightArrow{10.5} 
        \def\pad{0.15}          
        
        \def\yX{0}
        \pgfmathsetmacro{\yY}{\yX-\yGap}
        
        \filldraw[fill=gray!25,draw=black,thick] (\xL+4.5,\yX-0.5*\barH) rectangle ++(2.5,\barH);
        \filldraw[fill=gray!25,draw=black,thick] (\xL+4.5,\yY-0.5*\barH) rectangle ++(2.5,\barH);
        \draw[thick] (\xL,\yX-0.5*\barH) rectangle (\xR,\yX+0.5*\barH);
        \draw[thick] (\xL,\yY-0.5*\barH) rectangle (\xR,\yY+0.5*\barH);
        
        \draw (\xL+4.5,\yX-0.5*\barH) -- (\xL+4.5,\yX+0.5*\barH);
        \draw (\xL+7,\yX-0.5*\barH) -- (\xL+7,\yX+0.5*\barH);
        \node[anchor=west, scale = 0.6] at (\xR-0.7,\yX) {$T$};
        \node[anchor=west, scale = 0.6] at (\xL+5.7,\yX) {$t+\Delta$};
        \node[anchor=east, scale = 0.6] at (\xL+0.6,\yX) {$1$};
        \node[anchor=east, scale = 0.6] at (\xL+1,\yX) {$2$};
        \node[anchor=east, scale = 0.6] at (\xL+5,\yX) {$t$};
        
        \draw (\xL+4.5,\yY-0.5*\barH) -- (\xL+4.5,\yY+0.5*\barH);
        \draw (\xL+7,\yY-0.5*\barH) -- (\xL+7,\yY+0.5*\barH);
        \node[anchor=west, scale = 0.6] at (\xR-0.7,\yY) {$T$};
        \node[anchor=west, scale = 0.6] at (\xL+5.7,\yY) {$t+\Delta$};
        \node[anchor=east, scale = 0.6] at (\xL+0.6,\yY) {$1$};
        \node[anchor=east, scale = 0.6] at (\xL+1,\yY) {$2$};
        \node[anchor=east, scale = 0.6] at (\xL+5,\yY) {$t$};
        
        \node[anchor=east] at (\xL-0.5,\yX) {$\bX$};
        \node[anchor=east] at (\xL-0.5,\yY) {$\bY$};
        
        \draw[<->, thin] (\xL+5.5,\yX-0.5*\barH-\pad) -- (\xL+5.5,\yY+0.5*\barH+\pad);
        \node at ({0.2*(\xLeftArrow+\xRightArrow)}, {\yX}) {$\cdots$};
        \node at ({0.2*(\xLeftArrow+\xRightArrow)}, {\yY}) {$\cdots$};
        \node at ({0.8*(\xLeftArrow+\xRightArrow)}, {\yX}) {$\cdots$};
        \node at ({0.8*(\xLeftArrow+\xRightArrow)}, {\yY}) {$\cdots$};
        \end{tikzpicture}
    \end{minipage}
    \hfill
    \begin{minipage}[t]{0.48\textwidth}
        \centering
        \begin{tikzpicture}[>=Latex, thick, scale=0.55]
        \def\xL{0}       
        \def\xR{12}      
        \def\barH{0.6}   
        \def\yGap{2.2}   
        \def\xLeftArrow{1.0}   
        \def\xRightArrow{10.5} 
        \def\pad{0.15}          
        
        \def\yX{0}
        \pgfmathsetmacro{\yY}{\yX-\yGap}
        \filldraw[fill=gray!25,draw=black,thick] (\xL+6,\yX-0.5*\barH) rectangle ++(2.5,\barH);
        \filldraw[fill=gray!25,draw=black,thick] (\xL+2.8,\yY-0.5*\barH) rectangle ++(2.5,\barH);
        \draw[thick] (\xL,\yX-0.5*\barH) rectangle (\xR,\yX+0.5*\barH);
        \draw[thick] (\xL,\yY-0.5*\barH) rectangle (\xR,\yY+0.5*\barH);
        
        \draw (\xL+6,\yX-0.5*\barH) -- (\xL+6,\yX+0.5*\barH);
        \draw (\xL+8.5,\yX-0.5*\barH) -- (\xL+8.5,\yX+0.5*\barH);
        \node[anchor=west, scale = 0.6] at (\xR-0.7,\yX) {$T$};
        \node[anchor=west, scale = 0.6] at (\xL+7.2,\yX) {$i+\Delta$};
        \node[anchor=east, scale = 0.6] at (\xL+0.6,\yX) {$1$};
        \node[anchor=east, scale = 0.6] at (\xL+1,\yX) {$2$};
        \node[anchor=east, scale = 0.6] at (\xL+6.5,\yX) {$i$};
        
        \draw (\xL+2.8,\yY-0.5*\barH) -- (\xL+2.8,\yY+0.5*\barH);
        \draw (\xL+5.3,\yY-0.5*\barH) -- (\xL+5.3,\yY+0.5*\barH);
        \node[anchor=west, scale = 0.6] at (\xR-0.7,\yY) {$T$};
        \node[anchor=west, scale = 0.6] at (\xL+4,\yY) {$j+\Delta$};
        \node[anchor=east, scale = 0.6] at (\xL+0.6,\yY) {$1$};
        \node[anchor=east, scale = 0.6] at (\xL+1,\yY) {$2$};
        \node[anchor=east, scale = 0.6] at (\xL+3.3,\yY) {$j$};
        
        \node[anchor=east] at (\xL-0.5,\yX) {$\bX$};
        \node[anchor=east] at (\xL-0.5,\yY) {$\bY$};
        
        \draw[<->, thin] (\xL+7,\yX-0.5*\barH-\pad) -- (\xL+4,\yY+0.5*\barH+\pad);
        \node at ({0.3*(\xLeftArrow+\xRightArrow)}, {\yX}) {$\cdots$};
        \node at ({0.15*(\xLeftArrow+\xRightArrow)}, {\yY}) {$\cdots$};
        \node at ({0.9*(\xLeftArrow+\xRightArrow)}, {\yX}) {$\cdots$};
        \node at ({0.7*(\xLeftArrow+\xRightArrow)}, {\yY}) {$\cdots$};
        \end{tikzpicture}
    \end{minipage}
    \caption{Left panel: $\Psi_t$ is calculated using the temporally aligned time window $\{t, \cdots, t+ \Delta\}$, as illustrated by the shaded segments between $t$ and $t+ \Delta$. Right panel: $\Psi_{i,j}$ is calculated with the time-shifted windows, which are illustrated by the shaded segment between $i$ and $i+\Delta$ in $\bX$, and the shaded segment between $j$ and $j+\Delta$ in $\bY$.}
    \label{fig:detect_event_timeshift}
\end{figure}

To determine whether values as large as $\Psi_t$ could plausibly occur by random chance, we now use time-shifting to build an approximate null distribution against which we compare this observed value. For each $i,j$, we define
\[\Psi_{i,j} :=\psi\big( (X_i,\dots,X_{i+\Delta}), (Y_j,\dots,Y_{j+\Delta})\big).\]
See Figure~\ref{fig:detect_event_timeshift} for an illustration of the true observed evidence $\Psi_t$, and its time-shifted counterpart $\Psi_{i,j}$. 
Effectively, we are ``shifting'' the $X$ and $Y$ data streams relative to each other---while $\Psi_t$ measures agreement between $\bX$ and $\bY$ at (or around) time $t$, this new value $\Psi_{i,j}$ instead compares $\bX$ at time $i$ with $\bY$ at time $j$, which are not actually simultaneous measurements (when we choose $i\neq j$). The intuition is this: if $\Psi_t$ is large simply due to random chance (e.g., for LIGO, a vehicle was passing near the $X$ detector at time $t$, and by random chance, another vehicle was also passing near the $Y$ detector at time $t$), then it should be equally likely that $\Psi_{i,j}$ might be large for the same type of coincidental reason. On the other hand, a real event (such as a detected gravitational wave) at time $t$ can make $\Psi_t$ large, but a typical value of $\Psi_{i,j}$ (for $i\neq j$) would not exhibit the same behavior. Therefore, the $\Psi_{i,j}$ values form a ``control group'', to illustrate the values of $\Psi_t$ that might occur by random chance, without being caused by a true event detection.

With these definitions in place, we are now ready to compute a p-value for determining the significance of the observed evidence, $\Psi_t$:
\begin{equation}\label{eqn:pval_event}p_t = \frac{1}{(T-\Delta)^2}\sum_{i,j=1}^{T-\Delta}\One{\Psi_{i,j} \geq \Psi_t}.\end{equation}
This measures the proportion of the ``control group'' (the $\Psi_{i,j}$'s) that are greater than or equal to the actual observed evidence, $\Psi_t$. 
Intuitively, we would therefore hope the time-shifted p-value $p_t$ should indeed be a valid p-value, at least approximately. Our key question, then, is the following: for any rejection threshold $\alpha\in[0,1]$,
\begin{quote}
    Under what assumptions can we ensure that, in the absence of any true event, the probability of a false positive $\PP{p_t\leq \alpha}$ is bounded (ideally, approximately bounded by $\alpha$)?
\end{quote}
We may also be interested in multiple testing scenarios. If we are interested in searching for potential events across multiple time points, $t\in\cT$ (where $\cT\subseteq[T-\Delta]:=\{1,\dots,T-\Delta\}$ may be any subset of indices), we can use a standard Bonferroni correction to correct for multiple comparisons: we would declare an event at time $t\in\cT$ if $p_t\leq \alpha/|\cT|$, with the aim of bounding the total probability of any false discoveries by $\alpha$. We may then ask: for any rejection threshold $\alpha\in[0,1]$,
\begin{quote}
    Under what assumptions can we ensure that, in the absence of any true event, the probability of a false positive $\PP{\min_{t\in\cT}p_t\leq \alpha/|\cT|}$ is bounded (ideally, approximately bounded by $\alpha$)?
\end{quote}

\subsection{Theoretical guarantees}

By construction, the time-shifted p-value $p_t$ expresses the statistical significance of the observed evidence $\Psi_t$, as compared to the time-shifted values $\Psi_{i,j}$ that act as an approximation to the null distribution. However, can we be confident that this approximation indeed provides a valid p-value---particularly since, for time series data, we might have strong dependence between measurements within each data stream?

We next present theoretical guarantees verifying that, with only weak assumptions, the time-shifted p-value $p_t$ offers control of the chance of a false positive.

\subsubsection{A guarantee under stationarity}
Our first result places the fewest assumptions on the data: we assume only that each detector's data stream is stationary, i.e., its distribution does not change over time.

\begin{definition}\label{def:stationary}
    A time series $\bX = (X_1,\dots,X_T)$ is \emph{stationary} if, for every $1\leq s \leq T-1$, 
    \[(X_1,\dots,X_i) \eqd (X_{1+s},\dots,X_{i+s})\]
    holds for each $i=1,\dots,T-s$, where $\eqd$ denotes equality in distribution. 
\end{definition}
\noindent That is, stationarity of a time series $\bX$ implies that the joint distribution of $(X_1,\dots,X_i)$ (the first $i$ terms in the time series) is the same as the joint distribution of $(X_{1+s},\dots,X_{i+s})$ (which also consists of $i$ consecutive data values, but is shifted by $s$ time points).

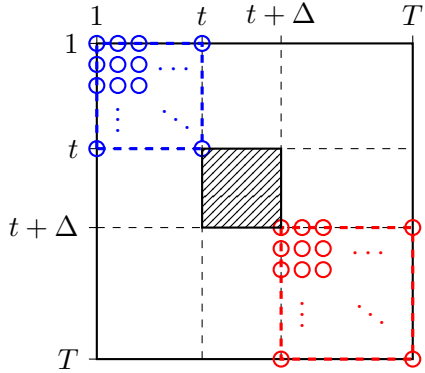
\begin{figure}
    \centering
        \begin{minipage}{0.48\textwidth}
        \centering
        \begin{tikzpicture}[scale=0.7,
    bluept/.style={draw=blue, circle, inner sep=2pt, thick},
    redpt/.style ={draw=red,  circle, inner sep=2pt, thick}
]

\draw[thick] (0,0) rectangle (6,6);

\foreach \x/\lab in {0/1,2/t,3.5/{t+\Delta},6/T}{
  \draw (\x,6) -- ++(0,0.15) node[above] {$\lab$};
}
\foreach \y/\lab in {0/T,2.5/{t+\Delta},4/t,6/1}{
  \draw (0,\y) -- ++(-0.15,0) node[left] {$\lab$};
}

\draw[dashed] (2,0) -- (2,6);
\draw[dashed] (3.5,0) -- (3.5,6);
\draw[dashed] (0,2.5) -- (6,2.5);
\draw[dashed] (0,4) -- (6,4);

\draw[blue,dashed, very thick] (0,4) rectangle (2,6);

\node[bluept] at (0,6){};
\node[bluept] at (0,5.6){};
\node[bluept] at (0,5.2){};
\node[bluept] at (0,4){};
\node[bluept] at (2,4){};
\node[bluept] at (2,6){};
\node[bluept] at (0.4, 6){};
\node[bluept] at (0.8, 6){};
\node[bluept] at (0.4, 5.6){};
\node[bluept] at (0.8, 5.6){};
\node[bluept] at (0.4, 5.2){};
\node[bluept] at (0.8, 5.2){};
\node at (0.4, 4.7) {\textcolor{blue}{$\vdots$}};
\node at (1.5, 5.5) {\textcolor{blue}{$\cdots$}};
\node at (1.5, 4.7) {\textcolor{blue}{$\ddots$}};

\draw[red,dashed, very thick] (3.5,2.5) rectangle (6,0);

\node[redpt] at (3.5,2.5) {};
\node[redpt] at (3.5,2.1) {};
\node[redpt] at (3.5,1.7) {};
\node[redpt] at (3.9,2.5) {};
\node[redpt] at (4.3,2.5) {};
\node[redpt] at (3.9,2.1) {};
\node[redpt] at (4.3,2.1) {};
\node[redpt] at (3.9,1.7) {};
\node[redpt] at (4.3,1.7) {};
\node[redpt] at (3.5,0) {};
\node[redpt] at (6,0) {};
\node[redpt] at (6,2.5) {};
\node at (3.9, 1) {\textcolor{red}{$\vdots$}};
\node at (5.2, 2) {\textcolor{red}{$\cdots$}};
\node at (5.2, 1) {\textcolor{red}{$\ddots$}};

\fill[gray!40,pattern=north east lines]
     (2,4) rectangle (3.5,2.5);
\draw[thick] (2,4) rectangle (3.5,2.5);

\end{tikzpicture}
    \end{minipage}
    \caption{An illustration of $\mathsf{mar}(t)$ as defined in~\eqref{eqn:define_mar(t)}. When computing $\mathsf{mar}(t)$, the numerator $\min\{t,T-\Delta+1-t\}^2$ is equal to minimum of the number of blue points, and the number of red points, shown in this figure. In contrast, the denominator $(T-\Delta)^2$ is approximately equal to the total number of points in the entire square.}
    \label{fig:illustrate_mar}
\end{figure}

The following result, Theorem~\ref{thm1}, offers a guarantee on the validity of the time-shifted p-value. Its guarantee is most useful for settings where the event time $t$ is near the middle of the range (i.e., $t\approx T/2$). Intuitively, this is important because it ensures that we have ample data both before and after time $t$ to use in our ``control group''. Before stating the theorem, we first define quantities that 
measure the size of the ``margin'' between the time range of interest (i.e., $\{t,\dots,t+\Delta\}$), and the endpoints of the full time range of the data (i.e., times $1$ and $T$):
let  \begin{equation}\label{eqn:define_mar(t)}\mathsf{mar}(t) =  \left(\frac{\min\{t , T-\Delta+1 -t\}}{T-\Delta}\right)^2.\end{equation}
In particular, if $t\approx T/2$ and $\Delta$ is small, we will have $\mathsf{mar}(t)\approx 1/4$. See Figure~\ref{fig:illustrate_mar} for an illustration.

\begin{theorem}\label{thm1}
    Assume $\bX = (X_1,\dots,X_T)$ and $\bY = (Y_1,\dots,Y_T)$ are each stationary time series. Then, if $\bX\independent\bY$, the time-shifted p-value $p_t$ defined in~\eqref{eqn:pval_event} satisfies
    \[\PP{p_t\leq \alpha} \leq \frac{\alpha}{\mathsf{mar}(t)},\]
    for any threshold $\alpha\in[0,1]$, and any choice of the function $\psi$. Consequently, for any (nonempty) subset of indices $\cT\subseteq[T-\Delta]$, it also holds that
    \[\PP{\min_{t\in\cT}p_t\leq \alpha/|\cT|} \leq \frac{\alpha}{\mathsf{mar}(\cT)},\]
    where $\mathsf{mar}(\cT) = \frac{|\cT|}{\sum_{t\in\cT} 1/\mathsf{mar}(t)}$ is the harmonic mean of the values $\{\mathsf{mar}(t)\}_{t\in\cT}$.
\end{theorem}

Note that we may have arbitrarily strong dependence within each stationary time series---for example, high autocorrelation between measurements $X_i$ and $X_j$ (particularly if the time points $i$ and $j$ are nearby in time). This means that the above theorem provides a guarantee of false positive control without placing any limitations on how much dependence there might be within each detector's time series. However, the upper bound is meaningful only if we are testing for an event at a time $t$ that is not too close to the endpoints of the time series---for instance, if $t\approx T/2$ then the inflation factor $\frac{1}{\mathsf{mar}(t)}$ is approximately $4$, but if $t\approx T$ then it may be arbitrarily large. In other words, the result above is meaningful only if we have ample data available both before and after time $t$, which may be true in certain applications but not in others.

In Appendix~\ref{app:thm1_factor}, we verify that this inflation, where the probability of a false positive may be larger than the nominal level $\alpha$, cannot be avoided without further assumptions. In particular, when $\alpha$ is small, the inflation factor of this upper bound is tight for any $t$: for example, if testing at the midpoint ($t\approx T/2$), it is indeed possible to obtain a false positive probability $\approx 4\alpha$ for certain distributions on the data.

\subsubsection{A guarantee under stationarity and $\beta$-mixing}
In some domains, we are specifically interested in ``real-time'' detection---that is, our task is to determine whether a real event is being detected at the present time, $t=T$, where we have available data from past times $t=1,\dots,T$. For this setting, the result of Theorem~\ref{thm1} is not useful, since the ``margin'' $\mathsf{mar}(t)$ is not bounded away from zero. 
Consequently, in order to obtain guarantees that are relevant for real-time detection, we will need to add an additional assumption that bounds the extent to which each detector's time series may have long-range dependence.
\begin{definition}\label{def:mixing}
    A time series $\bX = (X_1,\dots,X_T)$ is \emph{$\beta$-mixing}, with coefficient $\beta(\tau)$ at time lag $\tau\in\{0,\dots,T-1\}$, if it holds that
    \[\max_{i \in[T-\tau-1]}\dtv \big( 
    (X_1,\dots,X_i,X_{i+\tau+1},\dots,X_T) \, , \, (X_1,\dots,X_i,X'_{i+\tau+1},\dots,X'_T)\big)\leq \beta(\tau),\]
    where $\bX'$ is an independent copy of $\bX$ (i.e., $\bX' = (X'_1,\dots,X'_T)$ is drawn, independently, from the same joint distribution as $\bX$).
\end{definition}
\noindent In this definition, $\dtv$ denotes the total variation distance between distributions---namely, for distributions $P,Q$, $\dtv(P,Q) = \sup_A |P(A) - Q(A)|$, where the supremum is taken over all events $A$. We can interpret this definition as follows: if a time series $\bX$ is $\beta$-mixing, and its coefficient $\beta(\tau)$ at lag $\tau$ is relatively small, this means that the subsequences $(X_1,\dots,X_i)$ and  $(X_{i+\tau+1},\dots,X_T)$ are nearly independent. Note that $\bX$ may have strong correlation between observations at nearby times, but nonetheless we might have $\beta(\tau)$ small once the lag $\tau$ is sufficiently large.

Our second result provides a guarantee with this stronger assumption.

\begin{theorem}\label{thm2}
    Assume $\bX = (X_1,\dots,X_T)$ and $\bY = (Y_1,\dots,Y_T)$ are each stationary time series, and are $\beta$-mixing with coefficients $\beta_X(\tau)$ and $\beta_Y(\tau)$, respectively. Then, if $\bX\independent\bY$, the time-shifted p-value $p_t$ defined in~\eqref{eqn:pval_event} satisfies
   \begin{multline*}\PP{\min_{t\in\cT}p_t\leq \frac{\alpha}{|\cT|}}\leq {}\\ \alpha + \min_{\tau\geq 0}\left\{\min\left\{\sqrt{\alpha \cdot \frac{8(\tau+\Delta)}{T-\Delta}}
    ,(|\cT|-1)\cdot\frac{2(\tau+\Delta)}{T-\Delta}\right\}+ \frac{4(\tau+\Delta)}{T-\Delta}+ 2\beta_X(\tau)+2\beta_Y(\tau)\right\},\end{multline*}
    for any threshold $\alpha\in[0,1]$, any (nonempty) subset of indices $\cT\subseteq[T-\Delta]$, and any choice of the function $\psi$ with time window $\Delta$. 
\end{theorem}
\noindent In other words, this result ensures approximate control of the probability of a false detection (i.e., the upper bound is not much larger than $\alpha$), as long as the $\beta$-mixing coefficients are not too large (for some $\tau\ll T$).
As a special case, if we are only testing a single $p_t$ (i.e., not a multiple testing scenario), then this theorem yields
 \[\PP{p_t\leq \alpha} \leq \alpha 
    + \frac{4(\tau+\Delta)}{T-\Delta}+ 2\beta_X(\tau)+2\beta_Y(\tau)\]
    by taking $\cT=\{t\}$.

\paragraph{Comparing the theorems.}
Importantly, in the result of Theorem~\ref{thm2}, the upper bound does not vary with $t$---it offers a meaningful bound even for the ``real-time'' detection problem, when we aim to verify that an event has occurred at time $t=T$. Of course, this comes at the cost of requiring the additional $\beta$-mixing assumption, which is required to ensure that any dependence within a single data stream must be relatively short-lived. In Theorem~\ref{thm1}, on the other hand, since we do not assume any mixing conditions, then we need to be able to observe data occurring both before time $t$ and after time $t$ (i.e., $t$ should lie in the middle of the time range), to be able to assess whether the evidence observed at time $t$ is significant.

\section{Coincidence detection via time-shifting: detecting synchronicity}\label{sec:method_synchronicity}

We next turn to the question of detecting an overall pattern of synchronicity, and present a time-shifting method for this problem.

\subsection{The method}
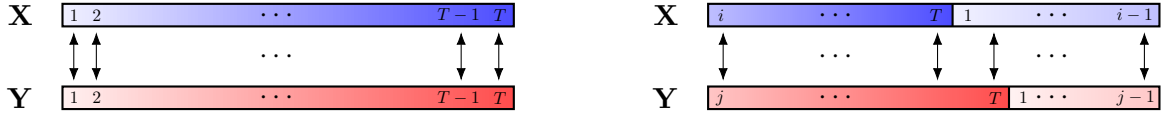
\begin{figure}
    \centering
        \begin{minipage}[t]{0.48\textwidth}
        \centering
        \begin{tikzpicture}[>=Latex, thick, scale=0.5]
        \def\xL{0}       
        \def\xR{12}      
        \def\barH{0.6}   
        \def\yGap{2.2}   
        \def\xLeftArrow{1.0}   
        \def\xRightArrow{10.5} 
        \def\pad{0.15}          
        
        \def\yX{0}
        \pgfmathsetmacro{\yY}{\yX-\yGap}
        \shade[left color=blue!5, right color=blue!70] (\xL,\yX-0.5*\barH) rectangle (\xR,\yX+0.5*\barH);
        \shade[left color=red!5, right color=red!70] (\xL,\yY-0.5*\barH) rectangle (\xR,\yY+0.5*\barH);
        \draw[thick] (\xL,\yX-0.5*\barH) rectangle (\xR,\yX+0.5*\barH);
        \draw[thick] (\xL,\yY-0.5*\barH) rectangle (\xR,\yY+0.5*\barH);
        
        \node[anchor=west, scale = 0.6] at (\xR-0.7,\yX) {$T$};
        \node[anchor=west, scale = 0.6] at (\xR-2.2,\yX) {$T-1$};
        \node[anchor=east, scale = 0.6] at (\xL+0.6,\yX) {$1$};
        \node[anchor=east, scale = 0.6] at (\xL+1.2,\yX) {$2$};
        
        \node[anchor=west, scale = 0.6] at (\xR-0.7,\yY) {$T$};
        \node[anchor=west, scale = 0.6] at (\xR-2.2,\yY) {$T-1$};
        \node[anchor=east, scale = 0.6] at (\xL+0.6,\yY) {$1$};
        \node[anchor=east, scale = 0.6] at (\xL+1.2,\yY) {$2$};
        
        \node[anchor=east] at (\xL-0.5,\yX) {$\bX$};
        \node[anchor=east] at (\xL-0.5,\yY) {$\bY$};
        
        \draw[<->, thin] (0.3,\yX-0.5*\barH-\pad) -- (0.3,\yY+0.5*\barH+\pad);
        \draw[<->, thin] (0.9,\yX-0.5*\barH-\pad) -- (0.9,\yY+0.5*\barH+\pad);
        \draw[<->, thin] (11.6,\yX-0.5*\barH-\pad) -- (11.6, \yY+0.5*\barH+\pad);
        \draw[<->, thin] (10.6,\yX-0.5*\barH-\pad) -- (10.6, \yY+0.5*\barH+\pad);
        \node at ({0.5*(\xLeftArrow+\xRightArrow)}, {0.5*(\yX+\yY)}) {$\cdots$};
        \node at ({0.5*(\xLeftArrow+\xRightArrow)}, {\yX}) {$\cdots$};
        \node at ({0.5*(\xLeftArrow+\xRightArrow)}, {\yY}) {$\cdots$};
        \end{tikzpicture}
    \end{minipage}
    \hfill
    \begin{minipage}[t]{0.48\textwidth}
        \centering
        \begin{tikzpicture}[>=Latex, thick, scale=0.5]
        \def\xL{0}       
        \def\xR{12}      
        \def\barH{0.6}   
        \def\yGap{2.2}   
        \def\xLeftArrow{1.0}   
        \def\xRightArrow{10.5} 
        \def\pad{0.15}          
        
        \def\yX{0}
        \pgfmathsetmacro{\yY}{\yX-\yGap}
        \shade[left color=blue!5, right color=blue!25] (\xL+6.5,\yX-0.5*\barH) rectangle ++(\xR-\xL-6.5,\barH);
        \shade[left color=red!5, right color=red!22] (\xR-4,\yY-0.5*\barH) rectangle ++(4,\barH);
        \shade[left color=blue!25, right color=blue!70] (\xL,\yX-0.5*\barH) rectangle ++(\xL+6.5,\barH);
        \shade[left color=red!22, right color=red!70] (\xL,\yY-0.5*\barH) rectangle ++(\xR-4,\barH);
        \draw[thick] (\xL,\yX-0.5*\barH) rectangle (\xR,\yX+0.5*\barH);
        \draw[thick] (\xL,\yY-0.5*\barH) rectangle (\xR,\yY+0.5*\barH);

        
        \draw (\xL+6.5,\yX-0.5*\barH) -- (\xL+6.5,\yX+0.5*\barH);
        \node[anchor=west, scale = 0.6] at (\xR-1.3,\yX) {$i-1$};
        \node[anchor=east, scale = 0.6] at (\xL+0.6,\yX) {$i$};
        \node[anchor=east, scale = 0.6] at (\xL+6.4,\yX) {$T$};
        \node[anchor=east, scale = 0.6] at (\xL+7.2,\yX) {$1$};
        
        \draw (\xR-4,\yY-0.5*\barH) -- (\xR-4,\yY+0.5*\barH);
        \node[anchor=west, scale = 0.6] at (\xR-1.3,\yY) {$j-1$};
        \node[anchor=east, scale = 0.6] at (\xL+0.6,\yY) {$j$};
        \node[anchor=west, scale = 0.6] at (\xR-4.7,\yY) {$T$};
        \node[anchor=west, scale = 0.6] at (\xR-3.9,\yY) {$1$};
        
        \node[anchor=east] at (\xL-0.5,\yX) {$\bX$};
        \node[anchor=east] at (\xL-0.5,\yY) {$\bY$};
        
        \draw[<->, thin] (0.4,\yX-0.5*\barH-\pad) -- (0.4,\yY+0.5*\barH+\pad);
        \draw[<->, thin] (\xL+6.1,\yX-0.5*\barH-\pad) -- (\xL+6.1,\yY+0.5*\barH+\pad);
        \draw[<->, thin] (11.6,\yX-0.5*\barH-\pad) -- (11.6, \yY+0.5*\barH+\pad);
        \draw[<->, thin] (\xR-4.4,\yX-0.5*\barH-\pad) -- (\xR-4.4, \yY+0.5*\barH+\pad);
        \node at ({0.3*(\xLeftArrow+\xRightArrow)}, {0.5*(\yX+\yY)}) {$\cdots$};
        \node at ({0.8*(\xLeftArrow+\xRightArrow)}, {0.5*(\yX+\yY)}) {$\cdots$};
        \node at ({0.3*(\xLeftArrow+\xRightArrow)}, {\yX}) {$\cdots$};
        \node at ({0.8*(\xLeftArrow+\xRightArrow)}, {\yX}) {$\cdots$};
        \node at ({0.8*(\xLeftArrow+\xRightArrow)}, {\yY}) {$\cdots$};
        \node at ({0.3*(\xLeftArrow+\xRightArrow)}, {\yY}) {$\cdots$};
        \end{tikzpicture}
    \end{minipage}
    \caption{Left panel: The test statistic $\Psi$ is calculated with temporally aligned data streams $\bX = (X_1,\dots,X_T)$ and $\bY = (Y_1,\dots,Y_T)$. Right panel: $\Psi_{i,j}$ is calculated using the ``wrap-around'' data streams $(X_i,\dots,X_T,X_1,\dots,X_{i-1})$ and $(Y_j,\dots,Y_T,Y_1,\dots,Y_{j-1})$.}
    \label{fig:detect_synchronicity_timeshift}
\end{figure}

For this setting, our goal is to measure whether the data streams $\bX$ and $\bY$ exhibit more synchronicity than would likely occur by random chance. In this section, then, we will now allow $\psi$---the function that provides a measure of evidence from the data---to depend on the entire duration of both time series, so that our test statistic is given by
\[\Psi = \psi\big((X_1,\dots,X_T),(Y_1,\dots,Y_T)\big).\]
For example, if our aim is to measure synchronicity between the spike trains of two neurons, we might choose the cross-correlation function, as in~\eqref{eq:Psi_crosscorr_neuro}. 
To determine the significance of the observed value $\Psi$, we will compare to time-shifted versions of the same test statistic: for each $i,j=1,\dots,T$, we define
\[\Psi_{i,j} = \psi\big((X_i,\dots,X_T,X_1,\dots,X_{i-1}),(Y_j,\dots,Y_T,Y_1,\dots,Y_{j-1})\big).\]
See Figure~\ref{fig:detect_synchronicity_timeshift} for an illustration of the true observed evidence $\Psi$, and its time-shifted counterpart $\Psi_{i,j}$. 

Here we can think of $(X_i,\dots,X_T,X_1,\dots,X_{i-1})$ as a ``wrap-around'' version of original data stream $\bX = (X_1,\dots,X_n)$, so that the $i$th entry $X_i$ is moved into the first position. Similarly, $\bY$ is wrapped around to produce $(Y_j,\dots,Y_T,Y_1,\dots,Y_{j-1})$. Note that, for $i\neq j$, these two data streams are misaligned---for instance, we are treating $X_i$ and $Y_j$ as though these data values occurred at the same time (i.e., at time $1$, since these values are in the first position of their respective wrap-around vectors), but in fact these measurements are from different times. As in the previous section, the intuition here is that any appearance of synchronicity that occurs simply due to random chance (e.g., the two neurons happen to fire at the same time) is equally likely to occur in the wrap-around version of the data. We then define the time-shifted p-value,
\begin{equation}\label{eqn:pval_overall}p = \frac{1}{T^2}\sum_{i=1}^T \sum_{j=1}^T\One{\Psi_{i,j} \geq \Psi}.\end{equation}

\subsection{Theoretical guarantee}
Next, we would like to provide a guarantee for the performance of the p-value $p$ defined in~\eqref{eqn:pval_overall}. However, since $\Psi$ now depends on the full data streams (rather than only a small portion of each data stream, as for detecting events in Section~\ref{sec:method_events}), this is not possible without further assumptions. 

To be able to establish a guarantee, we will need to introduce a small modification to our definition of the p-value, to make it slightly more conservative: we define
\begin{equation}\label{eqn:pval_overall_eps}p^{+\varepsilon} = \frac{1}{T^2}\sum_{i=1}^T \sum_{j=1}^T\One{\Psi_{i,j} + \varepsilon\geq \Psi},\end{equation}
where $\varepsilon>0$ is a small positive value. In other words, this value, $p^{+\varepsilon}$, can only be small if the measured evidence $\Psi$ is larger than most $\Psi_{i,j}$'s by a positive margin (where the choice of $\varepsilon$ determines this margin). We will also need to assume a stability property on the function $\psi$, to ensure that $\psi$ is not too sensitive to changing a small number of entries in the data streams. First we define some notation: for vectors $\bx=(x_1,\dots,x_T)$ and $\bx'=(x'_1,\dots,x'_T)$, we write $\bx\stackrel{\tau}{\asymp}\bx'$ if $\bx$ and $\bx'$ are equal everywhere except for a consecutive string of $\tau$ many entries, i.e., if for some $k\in\{1,\dots,T-\tau-1\}$,
\[x_i = x'_i \textnormal{ for all $i=1,\dots,k$ and for all $i=k+\tau+1,\dots,T$}\]
or, if for some $k\in\{0,\dots,\tau\}$,
\[x_i = x'_i \textnormal{ for all $i=k+1,\dots,k+T-\tau$}.\]
 
\begin{assumption}[Stability to perturbations of the data.]\label{asm:stability}
    The function $\psi$ satisfies
    the following bound:
   \begin{equation}\label{eqn:psi_k_l}\left|\psi(\bx,\by) - \psi(\bx',\by')\right| \leq \gamma(\tau)\textnormal{ for any $\bx,\bx',\by,\by'$ with $\bx\stackrel{\tau}{\asymp}\bx'$, $\by\stackrel{\tau}{\asymp}\by'$.}\end{equation}
\end{assumption}
\noindent That is, this coefficient bounds the amount by which $\psi(\bx,\by)$ can change, if we perturb data vectors $\bx$ and $\by$ by replacing $\tau$ many consecutive entries in each of the two data vectors. 

We are now ready to present our theoretical guarantee:
\begin{theorem}\label{thm3}
    Assume $\bX = (X_1,\dots,X_T)$ and $\bY = (Y_1,\dots,Y_T)$ are each stationary time series, and are $\beta$-mixing with coefficients $\beta_X(\tau)$ and $\beta_Y(\tau)$, respectively. Assume also that $\psi$ satisfies Assumption~\ref{asm:stability} with coefficients $\gamma(\tau)$. 
    Then, if $\bX\independent\bY$, 
    the time-shifted p-value $p^{+\varepsilon}$ defined in~\eqref{eqn:pval_overall_eps} satisfies
    \[\PP{p^{+\varepsilon}\leq \alpha} \leq \alpha +\min\left\{ 2\beta_X(\tau)+2\beta_Y(\tau) : \tau\geq 0, 4\gamma(\tau)\leq \varepsilon\right\},\]
    for any threshold $\alpha\in[0,1]$.
\end{theorem}
To interpret the upper bound, we note that $\gamma(0)=0$ for any function $\psi$, and therefore $4\gamma(\tau)\leq\varepsilon$ holds trivially for any $\varepsilon\geq 0$---but this may not offer a useful upper bound, since $\beta_X(0),\beta_Y(0)$ are typically large (i.e., consecutive measurements in the time series may have strong dependence). However, if $T$ is large and $\psi$ is not too sensitive to changing a small fraction of the input data, we would likely be able to find a value $\tau$ that is small enough to have $\gamma(\tau)\approx 0$ (i.e., $\tau \ll T$), but large enough to have $\beta_X(\tau),\beta_Y(\tau)\approx 0$ (i.e., $\tau \gg 1$), leading to a meaningful upper bound. 

In practice, when testing for synchronicity, we would typically choose a function $\psi$ that is fairly stable to small perturbations in the data, since we are aiming to measure the overall evidence for synchronicity across a long time range. For instance, the cross-correlation~\eqref{eq:Psi_crosscorr_neuro} measure typically changes only weakly if a small fraction of the values in $\bX$ and $\bY$ are changed. Consequently, we can expect $\gamma(\tau) \approx 0$ for moderate $\tau$ that ensures $\beta_X(\tau),\beta_Y(\tau)\approx 0$; therefore only a small, or even zero, inflation $\varepsilon$ would be needed.

Moreover, although the result of this theorem offers a guarantee that applies to a slightly more conservative $p^{+\varepsilon}$ defined in~\eqref{eqn:pval_overall_eps}, it is likely that the same type of false positive control holds for the unmodified p-value $p$ defined in \eqref{eqn:pval_overall} as well: in practical settings, if $\varepsilon$ is small, then we would expect to have $\One{\Psi_{i,j}\geq \Psi} = \One{\Psi_{i,j} + \varepsilon\geq \Psi}$, for most $i,j$, since this equality can only fail if the test statistic $\Psi$ is nearly tied with the ``control'' value $\Psi_{i,j}$.

\section{Empirical results}\label{sec:experiments}
In this section, we empirically evaluate the false detection rate of the time-shifting methods on both simulated data and real data.\footnote{Code for reproducing both simulated and real data experiments is available at \url{https://github.com/RuitingL/coincidence-detection}.}

\subsection{False positive control with simulated data}\label{sec:experiments_simulation}
We begin with a simulated data experiment to explore the false positive control properties of the time-shifted methods.

\paragraph{Data-generating process.}
For the simulated data experiment, we generate the data streams $\bX$ and $\bY$ independently, using the following stationary time series distribution.

Given parameters $\rho\in(0,1)$, $q \in (0,1)$ and $\sigma^2>0$, we first generate a mean vector $\mu\in\R^T$ via the following Markov chain: at each time $t\geq 1$,
\[\mu_{t+1} = \begin{cases} \mu_t + 1, & \textnormal{ with probability $q(1-\rho)$},\\
\mu_t - 1, & \textnormal{ with probability $(1-q)(1-\rho)$,}\\
0, & \textnormal{ with probability $\rho$.}\end{cases}\]
In other words, at each time step, with probability $1-\rho$ we follow a random walk with drift (if we choose $q>0.5$, then we will tend to see $\mu_t$ increase with $t$). But, with probability $\rho$, the process resets to $\mu_{t+1}=0$. The first time point, $\mu_1$, is drawn from the stationary distribution of the above Markov chain: we set $\mu_1 = 2B-N$, where $B\sim\textnormal{Binomial}(N,q)$, and where the random variable $N$ is drawn from the $\textnormal{Geometric}(\rho)$ distribution, with $\PP{N=k} = \rho(1-\rho)^k$ for integers $k\geq 0$. Finally, we return data obtained by adding Gaussian noise to $\mu$,
\[X_t \sim\mathcal{N}(\mu_t,\sigma^2),\]
with independent noise across times $t=1,\dots,T$.

\begin{figure}
    \centering
    \includegraphics[width=0.8\textwidth]{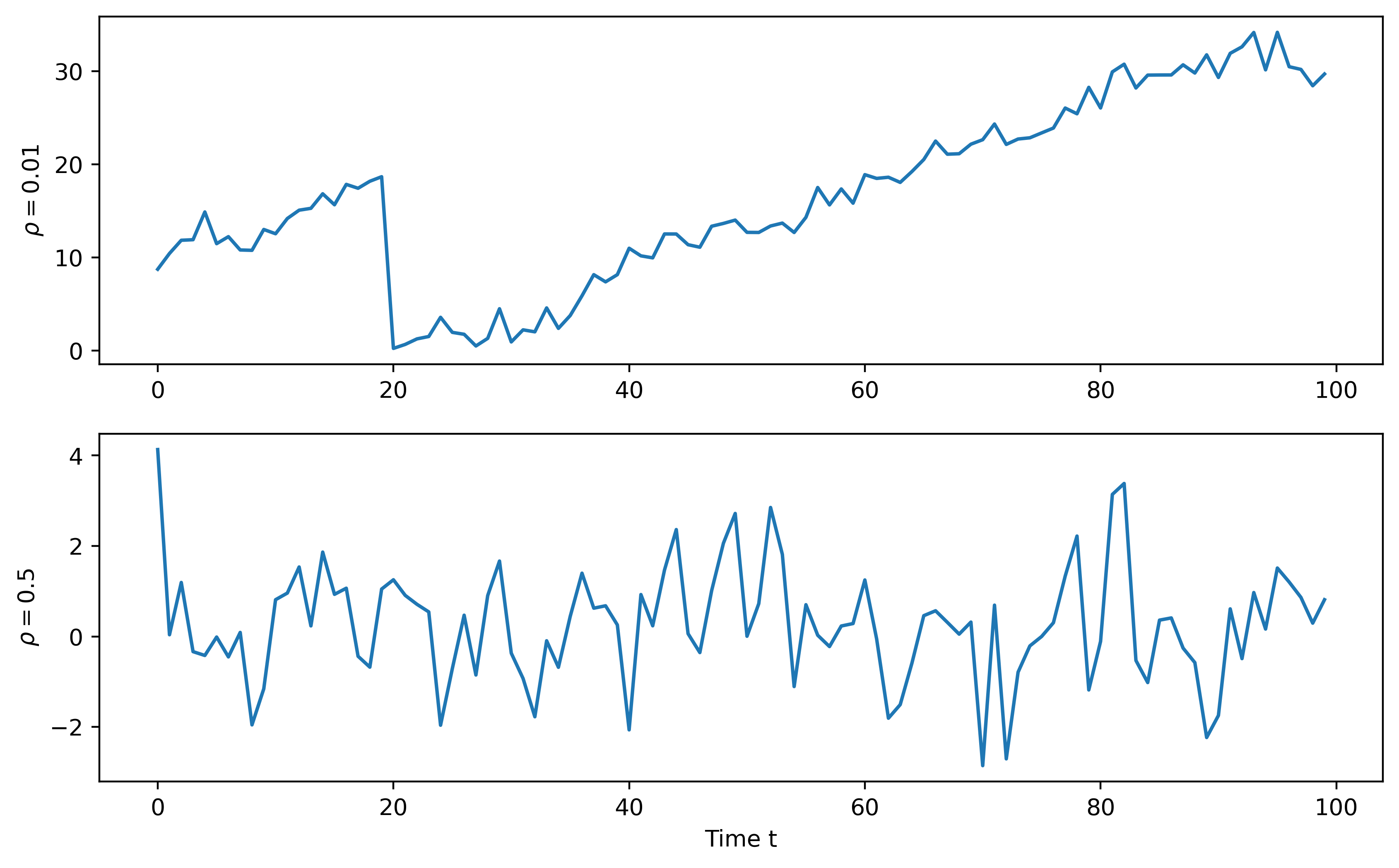}
    \caption{Illustration of a typical draw of the data for the simulated data experiment, with $\sigma=1$, $q=0.75$, and two different values of $\rho$. The top panel shows data generated with a small value $\rho=0.01$, corresponding to poor mixing of the time series---we see long periods of positive drift, indicating high temporal dependence.
    In contrast, the data in the bottom panel is generated with a larger value $\rho=0.5$, and the process appears to mix rapidly, with less dependence across time.
    }
    \label{fig:simulation_data_vis}
\end{figure}

In our simulations, we generate $\bX$ and $\bY$ independently from the above data-generating process, with time series length $T=100$. Each experiment is repeated for $N=10,000$ independent trials. For this experiment, we use parameters $\sigma=1$ and $q=0.75$, and vary the value of the parameter $\rho\in\{0.01,0.05,0.1,0.5\}$. The mixing rate of this process is controlled by this parameter $\rho$, with $\beta(\tau) \asymp (1-\rho)^{\tau+1}$; that is, small values $\rho\approx 0$ correspond to slower mixing.
Figure~\ref{fig:simulation_data_vis} shows two random draws of data streams from this process under two different $\rho$.

\begin{figure}
    \centering
    \includegraphics[width=\textwidth]{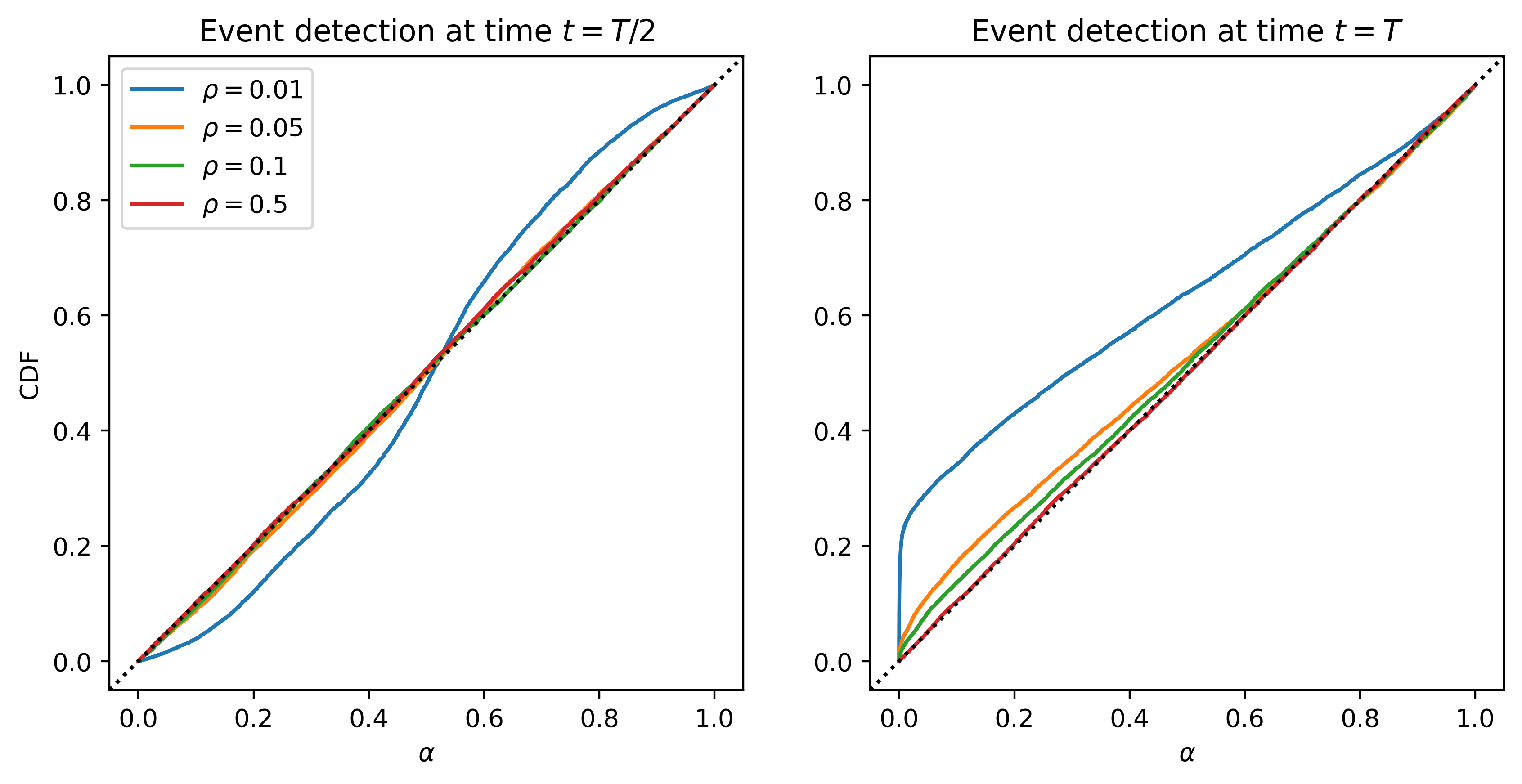}
    \caption{Results for detecting events with simulated data (see Section~\ref{sec:experiments_simulation} for details). For each value of $\rho$ used to generate the data, the figure shows the empirical cumulative distribution function (CDF) of the time-shifted p-value $p_t$ for detecting an event at time $t$, over $N=10000$ independent trials. The dotted diagonal line represents a uniform distribution.}
    \label{fig:simulation_event}
\end{figure}

\paragraph{Results for detecting events.}
For each trial, we test the time-shifted method for detecting events: following the notation of Section~\ref{sec:method_events}, we set $\Delta=0$, and use the function $\psi(x,y) = x+y$ as a measure of evidence. On each draw of the data, we compute the time-shifted p-value $p_t$ as in~\eqref{eqn:pval_event} for $t=T/2$, and again for $t=T$.

The results of this experiment are shown in Figure~\ref{fig:simulation_event}.
Since each experiment setting is repeated for $N=10,000$ independent draws of the subsampled data, we therefore observe $N$ values of the time-shifted p-value $p_t$ (for each choice $t=T/2$ and $t=T$). 
$\bX$ and $\bY$ are generated independently (i.e., the data satisfies the null hypothesis), and therefore we would want to see that the time-shifted p-value is uniformly distributed (or is more conservative), to avoid an inflated probability of a false positive detection. 

We can observe the following from our simulation results. For testing at the midpoint of the observation window in time ($t=T/2$), the false detection rate is well-controlled regardless of the mixing rate of the process, which is aligned with the finding in Theorem~\ref{thm1}. In contrast, when testing for an event at the endpoint $t=T$, the mixing properties play a key role: for poor mixing (when $\rho\approx 0$), it is not possible to provide any guarantees for the time-shifted p-value when $t$ is close to the endpoints of the time interval (since the results of Theorem~\ref{thm1} are not meaningful for this setting), which is illustrated by the steep CDF curve of $p_T$ for small $\rho$. On the other hand for larger values of $\rho$, the process mixes fairly rapidly, leading to approximate validity for detecting events even at time $t=T$, as suggested by the guarantee of Theorem~\ref{thm2}.

\paragraph{Results for detecting synchronicity.}
Next, for each draw of the data, we test the time-shifted method for detecting synchronicity: 
following the notation of Section~\ref{sec:method_synchronicity}, we use the function 
\[\psi(\bx,\by) = \frac{1}{w}\left(\Big(x_T + \sum_{t=1}^{w-1}(-1)^t x_t\Big)^2 + \Big(y_T + \sum_{t=1}^{w-1}(-1)^t y_t\Big)^2 \right) \]
as a measure of evidence, which has an additional parameter given by a positive integer $w\geq 1$. Note that this function is more stable (in the sense of Assumption~\ref{asm:stability}) for larger values of $w$, since each element in $\bx$ or $\by$ has less influence on the overall function value. We then compute the time-shifted p-value $p$ as in~\eqref{eqn:pval_overall}, and the modified p-value $p^{+\varepsilon}$ defined in~\eqref{eqn:pval_overall_eps} with inflation parameter $\varepsilon = 0.1$.

\begin{figure}
    \centering
    \includegraphics[width=\textwidth]{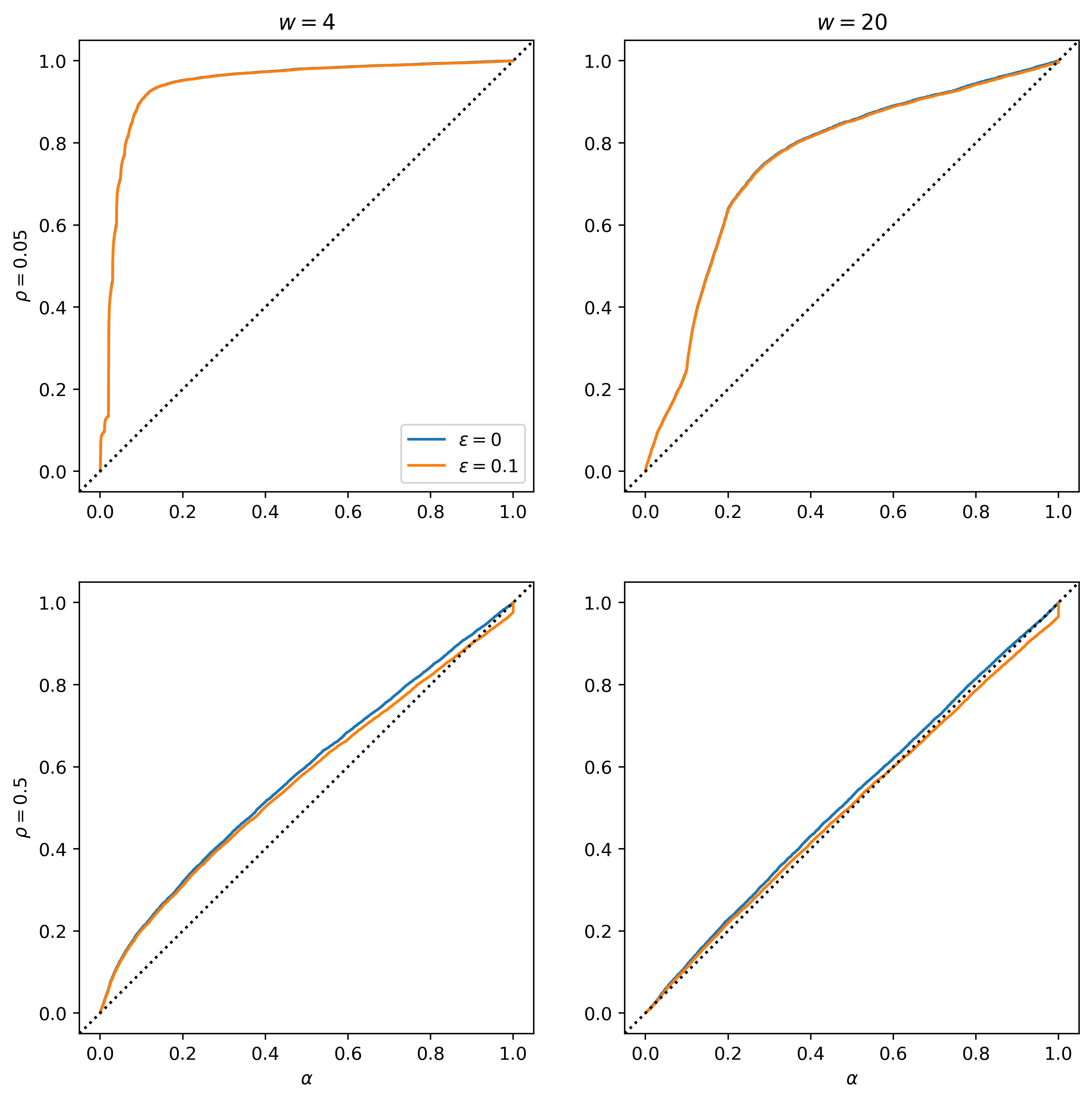}
    \caption{Results for detecting synchronicity with simulated data (see Section~\ref{sec:experiments_simulation} for details), with two different choices of the function $\psi$ that measures evidence of synchronicity. For each value of $\rho$ used to generate the data, the figure shows the empirical cumulative distribution function (CDF) of the time-shifted p-value $p$ for detecting synchronicity, over $N=10000$ independent trials, and the same for the inflated version of the p-value $p^{+\varepsilon}$ with $\varepsilon = 0.1$. The dotted diagonal line represents a uniform distribution.}
    \label{fig:simulation_sync}
\end{figure}

The results are displayed in Figure~\ref{fig:simulation_sync} (shown for $\rho = 0.05$ and $\rho=0.5$ only). First, as conjectured earlier, using a positive inflation parameter $\varepsilon=0.1$ does not noticeably alter the false positive probability; the uninflated p-value $p$ defined in~\eqref{eqn:pval_overall} has essentially the same performance as $p^{+\varepsilon}$, suggesting that the inflation may not be needed in practice. Next, we can observe that effective control of the false positive probability is achieved when the data distribution has better mixing (i.e., a larger value of $\rho$), and the function $\psi$ has better stability (i.e., a larger value of $w$), confirming the results of Theorem~\ref{thm3}.

\subsection{Real data}\label{sec:experiments_realdata}
We next present results of applying the time-shifting methods for both events detection and synchronicity detection on the stock price data streams. The dataset is sourced from the work of \citet{gibbs2021adaptive}.\footnote{The data can be downloaded from \url{https://github.com/isgibbs/DtACI/tree/main/MajorDailyReturns}.}

\begin{figure}
    \centering
    \includegraphics[width=\textwidth]{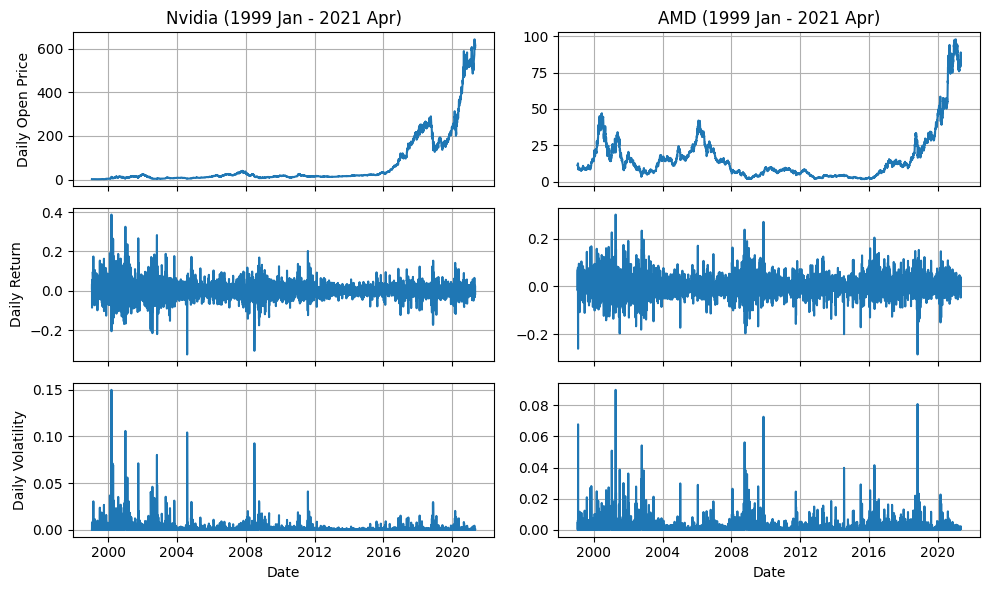}
    \caption{Stock price dataset (see Section~\ref{sec:experiments_realdata} for details).}
    \label{fig:Stockdata}
\end{figure}

\paragraph{Data.}
The data contains daily market open prices $P_t$ for two stocks, Nvidia and AMD, from January 1999 through April 2021. From these prices $P_t$, we also compute the daily simple return $R_t = \frac{P_t - P_{t-1}}{P_{t-1}}$ and volatility $V_t = R^2_t$. This data is measured over $T_{\textnormal{total}}=5604$ consecutive business days. The complete data is plotted in Figure~\ref{fig:Stockdata}. 

For our experiments, we will use data streams of length $T=100$. To generate the time series $\bX$ and $\bY$, we subsample the full data as follows:
\begin{itemize}
    \item Under the null (i.e., where the time series $\bX$ and $\bY$ are independent), we independently sample two indices $t_0,t_1\in\{0,\dots,T_{\textnormal{total}}-T\}$ and then define the time series as
    \[\bX = (P^{\textnormal{N}}_{t_0+1},\dots,P^{\textnormal{N}}_{t_0+T}), \quad \bY = (P^{\textnormal{A}}_{t_1+1},\dots,P^{\textnormal{A}}_{t_1+T}),\]
    where $P_t^{\textnormal{N}}$ and $P_t^{\textnormal{A}}$ denote the daily open price for Nvidia and for AMD, respectively, on day $t$ of the full dataset. Each time series contains $T$ consecutive time points and therefore may exhibit substantial temporal dependence; however, since the two time series are sampled at independent points in time (e.g., data in $\bX$ and $\bY$ may be measured in different years), it is reasonable to believe that $\bX$ and $\bY$ are independent.
    \item Under the alternative (i.e., where $\bX$ and $\bY$ may exhibit synchronicity or other shared dependence), we instead sample a single starting index $t_0\in\{0,\dots,T_{\textnormal{total}}-T\}$ at random, and define
    \[\bX = (P^{\textnormal{N}}_{t_0+1},\dots,P^{\textnormal{N}}_{t_0+T}), \quad \bY = (P^{\textnormal{A}}_{t_0+1},\dots,P^{\textnormal{A}}_{t_0+T}).\]
    Here, since $\bX$ and $\bY$ are measuring data from the same time window, we expect that $\bX$ and $\bY$ may exhibit some shared signal or synchronicity due to external events (e.g., overall market behavior).
\end{itemize}
We then repeat the above with returns data $R_t^{\textnormal{N}}$ and $R_t^{\textnormal{A}}$, and again with volatility data $V_t^{\textnormal{N}}$ and $V_t^{\textnormal{A}}$. In total, this leads to $6 = 3\times 2$ experiment settings, since we have three types of data (price, returns, volatility) each tested under the null and under the alternative. Each of these $6$ settings is repeated for $N=10000$ independent trials.

\begin{figure}
    \centering
    \includegraphics[width=\textwidth]{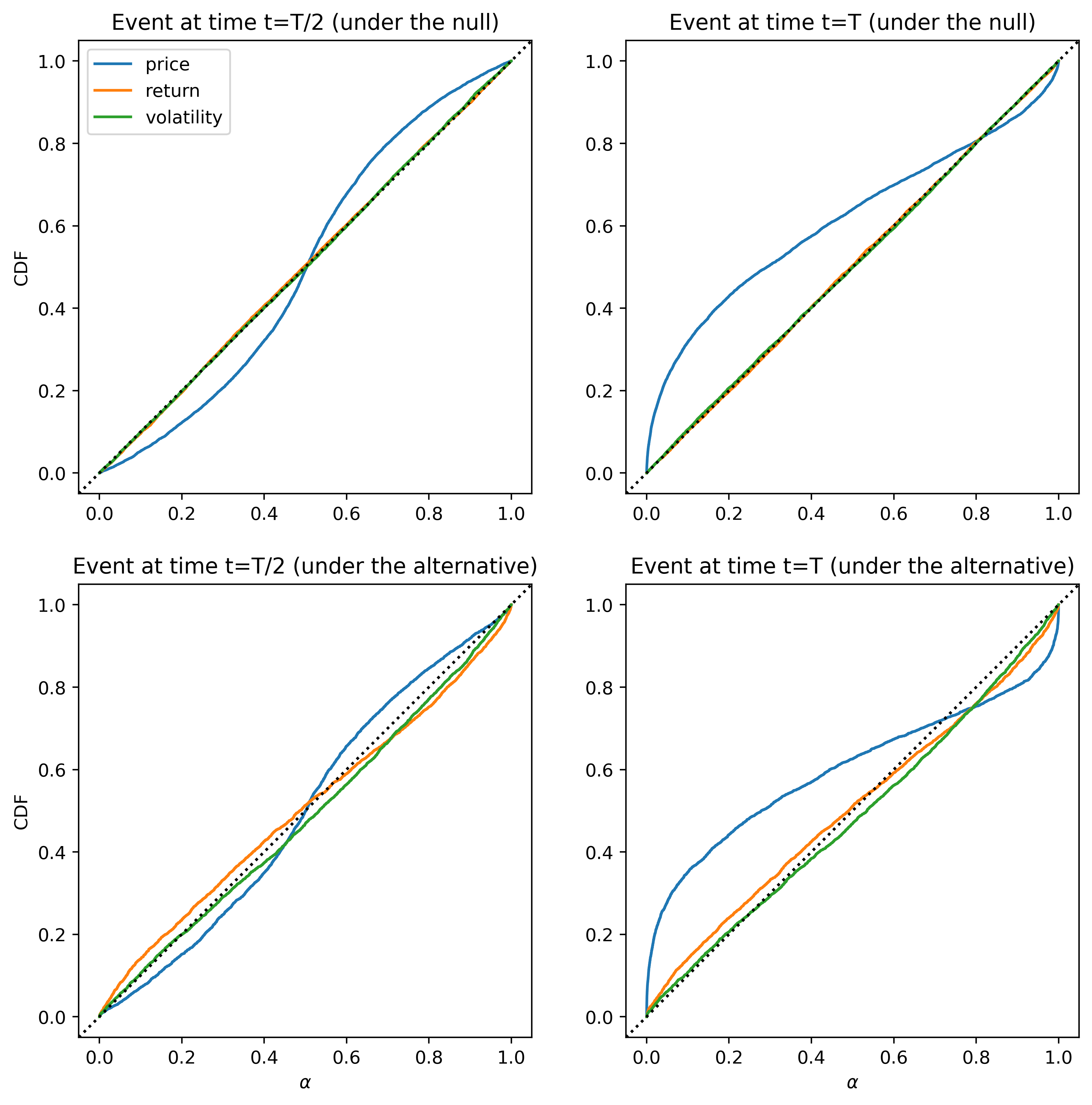}
    \caption{Results for detecting events with the stock price dataset (see Section~\ref{sec:experiments_realdata} for details). For each data set, the figure shows the empirical cumulative distribution function (CDF) of the time-shifted p-value $p_t$ for detecting an event at time $t$, over $N=10000$ independent trials. The dotted diagonal line represents a uniform distribution.}
    \label{fig:Stockdata_event}
\end{figure}

\paragraph{Results for detecting events.}
First, for each draw of the data, we test the time-shifted method for detecting events: following the notation of Section~\ref{sec:method_events}, we set $\Delta=0$, and use the function $\psi(x,y) = x+y$ as a measure of evidence. On each draw of the data, we compute the time-shifted p-value $p_t$ as in~\eqref{eqn:pval_event} for $t=T/2$, and again for $t=T$.

The results of this experiment are shown in Figure~\ref{fig:Stockdata_event}, which shows the CDF of the observed p-values across all trials, for each experiment.
Under the null hypothesis, ideally we would want to see that the time-shifted p-value is uniformly distributed (or is more conservative), to avoid an inflated probability of a false positive detection. Under the alternative, on the other hand, we would want small p-values to be more likely, to enable signal detection.

Testing at time $t=T/2$ shows the desired behavior under the null, with the CDF of the p-value lying near or below the diagonal line.  In contrast, when testing at time $t=T$, the price data shows extremely poor false positive control, while the returns data and volatility data continue to show good control of the false positive rate. This may be due to the fact that the price typically has stronger temporal dependence, as is evident in Figure~\ref{fig:Stockdata}, which corresponds to a worse mixing rate---that is, the price data fails to satisfy the mixing conditions required by Theorem~\ref{thm2} and thus only the results of Theorem~\ref{thm1} apply (which imply false positive control for $t=T/2$ but not for $t=T)$. In contrast, mixing conditions are plausibly approximately satisfied by the returns data and the volatility data, so Theorem~\ref{thm2} implies that false positive control holds for any $t$.

\paragraph{Results for detecting synchronicity.} Next, for each draw of the data, we test the time-shifted method for detecting synchronicity: 
following the notation of Section~\ref{sec:method_synchronicity}, we use the function $\psi(\bx,\by) = \textnormal{Corr}(\bx,\by)$ as a measure of evidence, where $\textnormal{Corr}$ is the sample correlation between two vectors. We then compute the time-shifted p-value $p$ as in~\eqref{eqn:pval_overall} (note that we do not use any inflation---that is, we do not use the modified p-value $p^{+\varepsilon}$ defined in~\eqref{eqn:pval_overall_eps}).

\begin{figure}
    \centering
    \includegraphics[width=\textwidth]{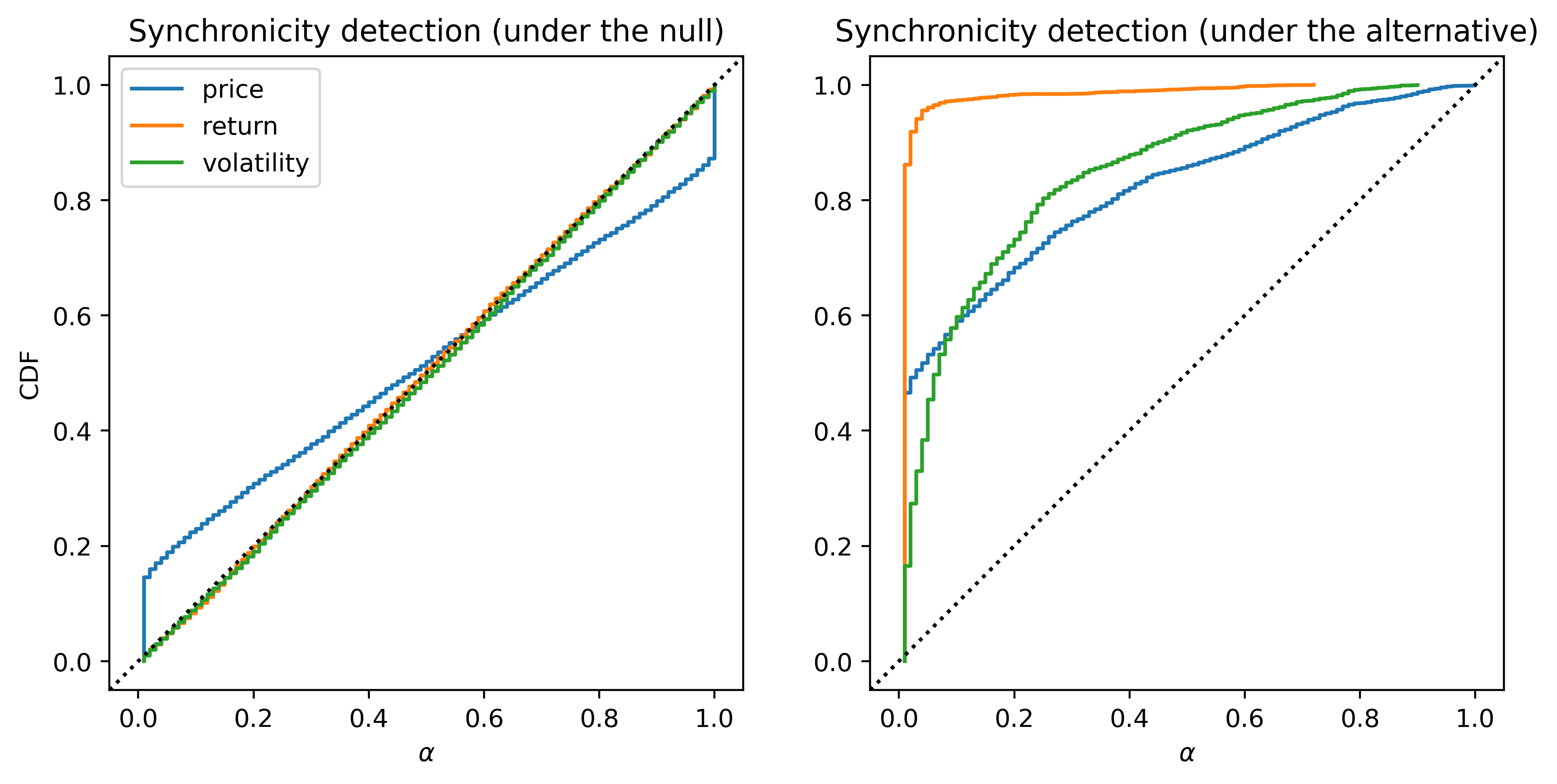}
    \caption{Results for detecting synchronicity with the stock price dataset (see Section~\ref{sec:experiments_realdata} for details). For each data set, the figure shows the empirical cumulative distribution function (CDF) of the time-shifted p-value $p$ for detecting synchronicity, over $N=10000$ independent trials. The dotted diagonal line represents a uniform distribution.}
    \label{fig:Stockdata_sync}
\end{figure}

The results of this experiment are shown in Figure~\ref{fig:Stockdata_sync}. As for detecting events, the ideal outcome is to have uniformly distributed p-values under the null, while under the alternative we would want a greater proportion of small p-values.

The results we see for synchronicity detection under the null are similar to our conclusions for detecting events. Namely, the price data (which shows high temporal dependence, and therefore, poor mixing) has extreme inflation of the false positive probability under the null, while the returns and volatility data show the desired behavior. This is consistent with the conclusions of Theorem~\ref{thm3} (note that we use the correlation function as the measure of evidence $\psi$, which is relatively stable as in Assumption~\ref{asm:stability}). Under the alternative, all three data types show strong signal, with a high probability of extremely low p-values, indicating high power to detect synchronicity.

\section{Discussion and related work}\label{sec:discussion}

In this work, we have examined the problem of coincidence detection in time series data, in settings where accurate models for the data distribution may not be available.
Our results show that time-shifting methods offer false positive control under only weak assumptions on the properties of the time series data. While the guaranteed upper bounds on the false positive probability do have some inflation above the target level $\alpha$, we emphasize that in practice we would often expect to see false positive probability $\approx \alpha$ (as is indeed the case in many of our simulation settings); we can interpret our results as providing reassurance that, even in the worst case, any inflation of the false positive probability cannot be too severe.
These findings support the existing use of time-shifting inference procedures across a range of scientific applications.

To conclude, we discuss some connections to prior work in the statistics literature, and describe some open questions 
and possible extensions of the methodology that would allow these results to be used for a broader range
of applications.

\subsection{Related work}
The methods explored in this work are closely related to many existing statistical methods, including classical methods such as permutation tests and bootstrapping, as well as more modern works on problems such as inference for time series and changepoint detection. 

\paragraph{Permutation tests of independence.}
Permutation testing is a classical tool that is frequently used to assess independence between multiple random variables: if the random vectors $\bX$ and $\bY$ consist of i.i.d.\ observations (rather than forming two time series, which implies that there may be dependence among the $X$ values and among the $Y$ values), we may choose random permutations $\sigma\in\mathcal{S}_T$ and $\sigma'\in\mathcal{S}_T$, and compare the true data vectors $\bX$ and $\bY$ (which may exhibit evidence of dependence) against their permuted copies $\bX_\sigma = (X_{\sigma(1)},\dots,X_{\sigma(T)})$ and $\bY_{\sigma'} = (Y_{\sigma'(1)},\dots,Y_{\sigma'(T)})$. Here, $\mathcal{S}_T$ is the symmetric group, i.e., the set of all $T!$ possible permutations of $[T]$. Permutation testing, both for the problem of testing independence between two random variables and for other inference problems, has a long history in statistic, dating back to classical work such as \citet{pitman1937significanceII,fisher1956mathematics}; see \cite{lehmann1986testing} for classical background on permutation testing.

In the context of astronomical time series data, \citet{Efron1992asimpletest, efron1994survival} construct a permutation test of independence for truncated data, and subsequently introduce a bootstrap approximation to the permutation p-values to reduce computational cost \citep{efron1999nonparametric}. In their work, however, the temporal structure of the sequences is not preserved, as all permutations within the time series are used to build the null distribution. In other words, using the notation of this present work, this line of the literature computes a p-value of the form
\[p = \frac{1}{(T!)^2} \sum_{\sigma,\sigma'\in\mathcal{S}_T} \One{\psi(\bX_\sigma, \bY_{\sigma'}) \ge \psi(\bX, \bY)},\]
or approximates this p-value by randomly sampling from the symmetric group $\mathcal{S}_T$. 

In contrast, the time-shifting method relies on p-values constructed as
\[p = \frac{1}{T^2} \sum_{i,j=1}^T \One{\psi(\bX_{[i]}, \bY_{[j]}) \ge \psi(\bX, \bY)},\]
which only uses ``control'' values of the form $\psi(\bX_{[i]},\bY_{[j]})$. Equivalently, defining the subset of circular permutations $\mathcal{S}_{T,\textnormal{circ}}\subseteq\mathcal{S}_T$, which consists of the $T$ permutations of the form
\[(1,2,\cdots,T) \mapsto (i,\cdots,T,1,\cdots,i-1), \quad i\in[T],\]
we are only considering ``control'' values $\psi(\bX_\sigma,\bY_{\sigma'})$ for $\sigma,\sigma'\in\mathcal{S}_{T,\textnormal{circ}}$. In contrast, classical permutation tests allow $\sigma,\sigma'$ to be \emph{any} permutation, may lead to ``control'' values $\psi(\bX_\sigma,\bY_{\sigma'})$ that do not reflect the null distribution of $\psi(\bX,\bY)$, since temporal dependence within the data streams is not preserved when using random permutations. 
In this sense, the time-shifting methods can be viewed as a special case of the general permutation test framework, designed to preserve temporal dependence. 

\paragraph{Inference via group actions, and applications to time series.} In our discussion above, we have seen that time-shifting methods can be viewed as a permutation test where we restrict to a subgroup of permutations, namely, $\mathcal{S}_{T,\textnormal{circ}}$. More generally, it is well known that permutation tests that average over any subgroup $G\subseteq\mathcal{S}_T$ are valid when the distribution of the test statistic $\psi(\bX,\bY)$ is invariant under the transformations in the subgroup $G$ \citep{hemerik2018exact}. Recent work by \citet{dobriban2025symmpi} extends this type of method to a more general framework, developing a general methodology for predictive inference when the data distribution exhibits group symmetries across arbitrary observation models.

\paragraph{Permutation tests and bootstrap tests for dependent data.}
The validity of permutation tests for dependent data has been studied for the problems of testing the dependence within a time series \citep{romano2022permutation} and testing whether $\bX$ is stationary \citep{romano2024permutation}. In these works, the asymptotic validity of the permutation test is derived under stationarity, $\alpha$-mixing and appropriate moment conditions on the sequence.
Moving to the context of testing the dependence of a pair of sequences $(\bX,\bY)$, \cite{romano2024least} investigate permutation tests for assessing the significance of (linear) regression coefficients in a time series setting, and show that the test is asymptotically valid under stationarity, $\alpha$-mixing and appropriate moment conditions.

In time-shifting methods for time series data, we use the subgroup of permutations $\mathcal{S}_{T,\textnormal{circ}}$ because permutations $\sigma\in\mathcal{S}_{T,\textnormal{circ}}$ approximately preserve the temporal dependence structure within a time series, and therefore, the distribution is approximately unchanged by applying a permutation $\sigma\in\mathcal{S}_{T,\textnormal{circ}}$. Many related procedures have been previously studied in the literature, for instance, permutations that partition the time series into blocks and then permute the blocks, as in the predictive inference work of \citet{chernozhukov2018exact}.
A closely related nonparametric method for generating null distributions is the bootstrap. This methodology was first introduced by \citet{efron1979bootstrap} for independent data, but was later extended to dependent data through block‐resampling techniques \citep{carlstein1986stationary, kunsch1989jackknife, politis1994stationary}, which preserve local dependence by resampling contiguous blocks of the sequence. Another variant is the circular bootstrap \citep{politis1991circular, politis1994stationary}, which uses circular permutations (i.e., permutations in $\mathcal{S}_{T,\textnormal{circ}}$) to construct ``wrap-around'' copies of the data as a preliminary step before performing block resampling, and is therefore closely related to the time-shifting methodology considered here.

\paragraph{Mixing conditions and valid inference for dependent data.}
Mixing conditions are a popular way to quantify weak dependence and are therefore widely used in establishing the validity of inference methods for dependent data. For example, in the permutation test setting \citep{romano2022permutation, romano2024least, romano2024permutation}; in the predictive inference setting \citep{chernozhukov2018exact, chernozhukov2021distributional, xu2023conformal, oliveira2024split, barber2025predictive}.
In particular, we note that the work of \citet{barber2025predictive} that examines the validity of conformal prediction when applied to time-dependent data, uses the $\beta$-mixing condition to show that the loss of coverage for conformal prediction is mild in time-series settings where the exchangeability does not hold. This is technically similar to our use of mixing conditions in establishing the validity of time-shifting methods, and some of our results use proof techniques similar to those developed in that work.

\paragraph{Changepoint detection.}
This paper focuses on coincidence detection based on two parallel observation sequences, where the two sequences corroborate each other to confirm the presence of a signal. In the broader context of signal detection, there is a large body of work on developing methods for detecting (distributional) change using a single observation sequence, notably in the area of changepoint detection. We highlight a few (recent) developments that are relevant because they are nonparametric and impose light distributional assumptions.
\citet{horvath1999testing} develop a test for changes in the mean of multivariate $m$-dependent stationary processes, which is asymptotically valid under a moment condition. 
For detecting general distributional changes, a recent line of work leverages conformal prediction for changepoint detection; see, for example, \citet{volkhonskiy2017inductive, vovk2021testing, nouretdinov2021conformal, dandapanthula2025conformal}.

\subsection{Extensions and applications}
The problems of detecting events, or detecting synchrony, arise in a broad range of applications such as astrophysics, neuroscience, and finance. As described in Section~\ref{sec:examples}, many of these fields already use time-shifting type methods, whether explicitly or implicitly, to determine the statistical significance of any detections. Our results offer some theoretical insight into the types of assumptions needed for such methods to be valid and provide reassurance that, even under very weak assumptions, the chance of a false positive can be controlled. However, there are many aspects of practical problems that are not yet captured in the framework considered in this present work. Here we briefly mention several important extensions of the methodology or of the setting, that would allow for broader applicability of these results:
\begin{itemize}
    \item \textbf{Detection in the streaming setting.} For the problem of detecting events, our results allow for error control in the setting of multiple testing, with Theorems~\ref{thm1} and~\ref{thm2} providing guarantees for $\min_t p_t$ where $p_t$ is the time-shifted p-value assessing the evidence for an event that occurred at time $t$. However, each of these values $p_t$ is computed using the entire data stream, over the time range $1,\dots,T$, and is therefore working in the ``batch'' setting: after observing data over a duration of $T$ time points, we then evaluate the evidence across all (or many) times $t\in\{1,\dots,T-\Delta\}$. In practice, we might instead want to perform real-time assessment at each time $t$, with $p_t$ computed using only the data observed prior to time $t$. Can the results of these theorems be extended to control the multiple testing problem in this streaming setting?
    \item \textbf{Real-time updates to the function $\psi$.} Continuing with the motivation of a streaming setting as described above, the function $\psi$ itself, which we use as a measure of evidence at time $t$, might be constructed in a data-dependent way. Consequently, in a streaming setting, this function $\psi$ might in fact vary with time---for instance, at time $t$, we evaluate the data $X_t,Y_t$ using a function $\psi$ that was trained or tuned based on data observed during times $1,\dots,t-1$. Are any guarantees of false positive control possible in this data-adaptive setting?
    \item \textbf{Deleting the diagonal from the control group.} Our p-value $p_t$ for detecting an event at time $t$~\eqref{eqn:pval_event} compares the test statistic $\Psi_t = \psi(X_t,Y_t)$ against ``control'' values $\Psi_{i,j} = \psi(X_i,Y_j)$, for all $i,j$ (here we are considering the case $\Delta=0$ for simplicity). However, if true events do occur throughout the time range, for any pair $i\approx j$, these ``control'' values may be larger than what we would expect simply by random chance because $X_i$ and $Y_j$ may both be affected by events occurring around time $i\approx j$. We might therefore choose to compute $p_t$ by comparing $\Psi$ only against values $\Psi_{i,j}$ for pairs $i,j$ sufficiently separated in time. Similarly, for detecting synchronicity, for $i\approx j$, a ``control'' value $\Psi_{i,j} = \psi((X_i,\dots,X_T,X_1,\dots,X_{i-1}),(Y_j,\dots,Y_T,Y_1,\dots,Y_{j-1}))$ may likely be larger than we would expect simply due to random chance, since the two time-shifted sequences are nearly aligned in time and may be exhibiting evidence of real synchronicity between the two data streams. Consequently, when computing $p$ as in~\eqref{eqn:pval_overall} (or $p^{+\varepsilon}$ as in~\eqref{eqn:pval_overall_eps}), we again might choose to exclude pairs $i,j$ with $i\approx j$. While these modifications of the method might be expected to have higher power (if true signals are present in the data), would these modified versions of the method offer similar false positive control guarantees?
    \item \textbf{Data-adaptive measures of stability.} For the problem of detecting synchronicity, the stability condition on $\psi$ stated in Assumption~\ref{asm:stability} requires that $\psi$ is not sensitive to small changes in \emph{any} inputs $\bx,\by$. While this is true for many practical choices of $\psi$, it is also likely that we might be able to relax this condition for broader applicability, by defining data-adaptive or average-case notions of stability: that is, would the same types of guarantees hold if we only require that $\psi$ is not sensitive to small changes in a ``typical'' draw of the data $\bX,\bY$?
\end{itemize}

\subsection*{Acknowledgements} 
We gratefully acknowledge the support of the NSF-Simons AI-Institute for the Sky (SkAI) via grants NSF AST-2421845 and Simons Foundation MPS-AI-00010513.
Additionally, R.L. and R.F.B. were partially supported by the National Science Foundation via grant DMS-2023109, R.F.B. was partially supported by the Office of Naval Research via grant N00014-24-1-2544. S.D. and D.E.H. were supported by the National Science Foundation via grants PHY-2110507 and PHY-2513312. D.E.H. was also supported by the Simons Collaboration on Black Holes and Strong Gravity through grant SFI-MPS-BH-00012593-07.
The authors thank Brent Doiron, Margaret Gardel, and Elizabeth Jerison for helpful discussions.

\clearpage

\appendix

\section{Proofs}
\subsection{Proof of Theorem~\ref{thm1}}
First, define
\[T' = \min\{t,T-\Delta+1-t\},\]
and fix any indices 
\[k,\ell\in\{1,\dots,T'\}.\]
Define a partial version of the p-value,
\[p_{t,k,\ell} = \frac{1}{T'{}^2}\sum_{i=t-k+1}^{t-k+T'}\sum_{j=t-\ell+1}^{t-\ell+T'}\One{\Psi_{i,j}\geq\Psi_t}.\]
Note that, by construction, it holds deterministically that
\[p_t = \frac{1}{(T-\Delta)^2}\sum_{i,j=1}^{T-\Delta}\One{\Psi_{i,j}\geq\Psi_t} \geq \frac{T'{}^2}{(T-\Delta)^2}p_{t,k,\ell}=\mathsf{mar}(t) \cdot  p_{t,k,\ell},\]
where the last step holds since $\frac{T'{}^2}{(T-\Delta)^2}  = \mathsf{mar}(t)$, by construction.
Therefore,
\[\PP{p_t\leq \alpha}\leq \PP{p_{t,k,\ell}\leq \frac{\alpha}{\mathsf{mar}(t)}}.\]
Next, note that $p_{t,k,\ell}$ depends on $\bX$ and $\bY$ only through the subvectors $(X_{t-k+1},\dots,X_{t-k+T'+\Delta})$ and $(Y_{t-\ell+1},\dots,Y_{t-\ell+T'+\Delta})$. Since $\bX$ and $\bY$ are each stationary and are independent, the distribution of $p_{t,k,\ell}$ is invariant to shifts of these time windows---and in particular,
\[p_{t,k,\ell} \eqd \frac{1}{T'{}^2}\sum_{i,j=1}^{T'}\One{\Psi_{i,j}\geq\Psi_{k,\ell}}=:p^*_{k,\ell}.\]
Thus, for any $k,\ell$,
\[\PP{p_t\leq \alpha}\leq \PP{p^*_{k,\ell}\leq \frac{\alpha}{\mathsf{mar}(t)}}.\]
But since this is true for any fixed $k,\ell$,
 the same is true if we sample indices $K,L$ randomly: drawing
\[K,L\iidsim\textnormal{Unif}\big([T']\big),\]
we have
\[\PP{p_t\leq \alpha}\leq \PP{p^*_{K,L}\leq \frac{\alpha}{\mathsf{mar}(t)}}.\]
And, $p^*_{K,L}$ is clearly a valid p-value, since it compares $\Psi_{K,L}$ against the empirical distribution of the values $(\Psi_{i,j})_{i,j}\in\R^{T'\times T'}$, when $\Psi_{K,L}$ is itself a uniform draw from this array. Thus
\[\PP{p^*_{K,L}\leq \frac{\alpha}{\mathsf{mar}(t)}}\leq \frac{\alpha}{\mathsf{mar}(t)},\]
which completes the proof for a single time $t$. Finally, for the setting of multiple comparisons, we simply apply a union bound,
 \[\PP{\min_{t\in\cT}p_t\leq \alpha/|\cT|} \leq \sum_{t\in\cT}\PP{p_t\leq \alpha/|\cT|}
 \leq   \sum_{t\in\cT} \frac{\alpha/|\cT|}{\mathsf{mar}(t)} = \frac{\alpha}{\mathsf{mar}(\cT)},\]
 where the second step holds by applying the first part of the theorem (with $\alpha/|\cT|$ in place of $\alpha$).

\subsection{Preliminaries for the proofs of Theorems~\ref{thm2} and~\ref{thm3}}
We begin by defining some notation for various subvectors or transformations of the data streams,
which will be useful for proving the theorems.

First, for any $\bx=(x_1,\dots,x_T)$, and any $i\in[T]$, we define 
the ``wrap-around'' data stream
\[\bx_{[i]} = (x_i,\dots,x_T,x_1,\dots,x_{i-1}).\] 
Next, for any $k\in[T]$, we define a subset of indices $I_{k,\tau}$, as follows: we define
\[I_{k,\tau} = \{1,\dots,k,k+\tau+1,\dots,T\}\textnormal{ if $1\leq k\leq T-\tau-1$},\]
where we note that the $\tau$ many indices immediately following time $k$ have been deleted, and we also define
\[I_{k,\tau} = \{k-T+\tau+1,\dots,k\}\textnormal{ if $T-\tau\leq k\leq T$},\]
where again $\tau$ many indices have been deleted, except that now these deletions occur at both the beginning and the end of the time range.
(Recalling the notation $\bx\stackrel{\tau}{\asymp}\bx'$ from Section~\ref{sec:method_synchronicity}, we can observe that 
$\bx\stackrel{\tau}{\asymp}\bx'$ holds if and only if $\bx_{I_{k,\tau}}=\bx'_{I_{k,\tau}}$ for some $k\in[T]$.)
The next lemma verifies a basic property of these sets. By an elementary counting argument, it holds that
\begin{equation}\label{eqn:counting}
    \textnormal{For any $(i,j)\in[T]$, }\sum_{k,\ell=1}^T \One{(i,j)\in I_{k,\tau}\times I_{\ell,\tau}} = (T-\tau)^2.
\end{equation}

We next compare the distributions of the data streams $\bX,\bY$ with the distributions of their ``wrap-around'' versions, as follows:
\begin{lemma}\label{lem:dtv_rotate}
    Let $\bX,\bY$ be stationary time series with $\beta$-mixing coefficients $\beta_X(\tau),\beta_Y(\tau)$, respectively. For any $i,j\in[T]$, it holds that
    \[\dtv\big(\bX_{I_{T-i+1,\tau}\cap I_{T,\tau}}, (\bX_{[i]})_{I_{T-i+1,\tau}\cap I_{T,\tau}}\big)\leq 2\beta_X(\tau)\]
    and similarly
    \[\dtv\big(\bY_{I_{T-j+1,\tau}\cap I_{T,\tau}}, (\bY_{[j]})_{I_{T-j+1,\tau}\cap I_{T,\tau}}\big)\leq 2\beta_Y(\tau).\]
\end{lemma}
Since $\bX\independent\bY$, consequently we have
\begin{multline}\label{eqn:dtv_for_thm3}\dtv\Big(\big(\bX_{I_{T-i+1,\tau}\cap I_{T,\tau}},\bY_{I_{T-j+1,\tau}\cap I_{T,\tau}}\big), \big((\bX_{[i]})_{I_{T-i+1,\tau}\cap I_{T,\tau}},(\bY_{[j]})_{I_{T-j+1,\tau}\cap I_{T,\tau}}\big)\Big)\leq 2\beta_X(\tau)+2\beta_Y(\tau),\end{multline}

Next, let $\cX,\cY$ denote the original data spaces, i.e., $(X_1,\dots,X_T)\in\cX^T$ and $(Y_1,\dots,Y_T)\in\cY^T$.
We also introduce a new symbol, $\varnothing$, which we can think of as representing ``missing data''.
For any $k\in[T]$, define
\[\bX^k = (X_1,\dots,X_k,\varnothing,\dots,\varnothing,X_{k+\tau+1},\dots,X_T),\]
and
\[\bX_k = (X_{T-k+1},\dots,X_T,\varnothing,\dots,\varnothing,X_1,\dots,X_{T-k-\tau}),\]
if $k\leq T-\tau-1$, or 
\[\bX^k = (\underbrace{\varnothing,\dots,\varnothing}_{\hspace{-.08in}k-T+\tau \textnormal{ elements}\hspace{-.08in}},X_{k-T+\tau+1},\dots,X_k,\underbrace{\varnothing,\dots,\varnothing}_{T-k \textnormal{ elements}}),\]
and
\[\bX_k = (\underbrace{\varnothing,\dots,\varnothing}_{\hspace{-.08in}k-T+\tau \textnormal{ elements}\hspace{-.08in}},X_{T-k+1},\dots,X_{2T-k-\tau},\underbrace{\varnothing,\dots,\varnothing}_{T-k \textnormal{ elements}}),\]
if $k\geq T-\tau$. Note that, in both cases, $(\bX^k)_i$ and $(\bX_k)_i$ take values in $\cX$ if $i\in I_{k,\tau}$,
and are equal to $\varnothing$ if $i\not\in I_{k,\tau}$.
Define $\bY^\ell$ and $\bY_\ell$ analogously for $\ell\in[T]$.
By definition of $\beta$-mixing, it holds that
\[\dtv(\bX^k,\bX_k)\leq 2\beta_X(\tau), \quad \dtv(\bY^k,\bY_k)\leq 2\beta_Y(\tau),\]
for each $k,\ell$. Since $\bX\independent \bY$, we therefore have
\begin{equation}\label{eqn:dtv_for_thm2}
\dtv\big((\bX^k,\bY^\ell),(\bX_k,\bY_\ell)\big)\leq 2\beta_X(\tau)+2\beta_Y(\tau).
\end{equation}

\subsection{Proof of Theorem~\ref{thm2}}
We will split into cases: first we will prove the result in the special case $\Delta=0$, and then
we will reduce the general case $\Delta\geq 1$ to this special case.

\subsubsection{Proof for the case $\Delta=0$}
For any $x\in\cX\cup\{\varnothing\}$, $y\in\cY\cup\{\varnothing\}$, define
\[\psi^*(x,y) = \begin{cases} \psi(x,y), & \textnormal{ if $x\in\cX$ and $y\in\cY$},\\ -\infty, &\textnormal{ otherwise},\end{cases}\]
i.e., $\psi^*(x,y)=-\infty$ if either $x=\varnothing$ or $y=\varnothing$. Define also
\[Q_a(\bx,\by) = \max\left\{q : \frac{1}{T^2}\sum_{i,j=1}^T\One{\psi^*(x_i,y_j)\geq q} \geq a\right\},\]
for any $a\in[0,1]$, which is a $(1-a)$-quantile of the values $(\psi^*(x_i,y_j))_{i,j\in[T]}$.
By construction, we can verify that for any $t\in\cT$,
\[p_t \leq a \ \Longleftrightarrow \ \psi(X_t,Y_t) \geq Q_a(\bX,\bY)  \ \Longleftrightarrow \ \psi^\ast(X_t,Y_t) \geq Q_a(\bX,\bY).\]
Define also
\[f(\bx,\by) = \max_{t\in\cT} \psi^*(x_t,y_t)  =\max_{t\in\cT} \psi(x_t,y_t),\]
so that we have
\[\min_{t\in\cT}p_t \leq a \ \Longleftrightarrow \ f(\bX,\bY) \geq Q_a(\bX,\bY), \]
for any $a\in[0,1]$. From this point on, then, our goal is to bound $\PP{f(\bX,\bY) \leq Q_{\alpha/|\cT|}(\bX,\bY)}$.

Now fix any $k,\ell\in[T]$. Since $\psi^*((\bX^k)_i,(\bY^\ell)_j)\leq \psi^*(X_i,Y_j)$ for all $i,j$ by construction, we therefore have $Q_a(\bX,\bY) \geq Q_a(\bX^k,\bY^\ell)$ and $f(\bX,\bY)\geq f(\bX^k,\bY^\ell)$. Define also a random variable $t^*$, which achieves the maximum, i.e.,
\[t^*\in\arg\max_{t\in\cT} \psi(X_t,Y_t).\]
If $(\bX^k)_{t^*} = X_{t^*}, (\bY^\ell)_{t^*}=Y_{t^*}$ (which holds whenever $t^*\in I_{k,\tau}\cap I_{\ell,\tau}$), we then have 
\[f(\bX,\bY) = \max_{t\in\cT} \psi(X_t,Y_t) = \psi(X_{t^*},Y_{t^*}) = \psi((\bX^k)_{t_*},(\bY^\ell)_{t^*}) \leq f(\bX^k,\bX^\ell),\]
i.e., we have shown that
\[t^*\in I_{k,\tau}\cap I_{\ell,\tau} \ \Longrightarrow \ f(\bX,\bY)=f(\bX^k,\bY^\ell).\]
Combining these calculations, and setting $a=\alpha/|\cT|$, we then have
\begin{align*}
    \PP{\min_{t\in\cT}p_t\leq \frac{\alpha}{|\cT|}} 
    & = \PP{f(\bX,\bY)\geq Q_{\alpha/|\cT|}(\bX,\bY)}\\
    &\leq\PP{f(\bX^k,\bY^\ell)\geq Q_{\alpha/|\cT|}(\bX^k,\bY^\ell)} + \PP{t^*\not\in  I_{k,\tau}\cap I_{\ell,\tau}}\\
    & \leq\PP{f(\bX_k,\bY_\ell)\geq Q_{\alpha/|\cT|}(\bX_k,\bY_\ell)} + \PP{t^*\not\in  I_{k,\tau}\cap I_{\ell,\tau}} + 2\beta_X(\tau)+2\beta_Y(\tau)\\
    &\leq \PP{f(\bX_k,\bY_\ell)\geq Q_{\alpha'/|\cT|}(\bX,\bY)} + \PP{t^*\not\in  I_{k,\tau}\cap I_{\ell,\tau}} \\
    &\hspace{1in}{} + 2\beta_X(\tau)+2\beta_Y(\tau)+ \PP{Q_{\alpha'/|\cT|}(\bX,\bY) > Q_{\alpha/|\cT|}(\bX_k,\bY_\ell)},
\end{align*}
where $\alpha'>\alpha$ will be specified later, and where the next-top-last step applies the total variation bound computed in~\eqref{eqn:dtv_for_thm2}. 
Since these calculations hold for any $k,\ell$, we can then take an average, to obtain:
\begin{multline*}
    \PP{\min_{t\in\cT}p_t\leq \frac{\alpha}{|\cT|}}
    \leq \frac{1}{T^2}\sum_{k,\ell=1}^T \PP{f(\bX_k,\bY_\ell)\geq Q_{\alpha'/|\cT|}(\bX,\bY)}
    + \frac{1}{T^2}\sum_{k,\ell=1}^T \PP{t^*\not\in  I_{k,\tau}\cap I_{\ell,\tau}}\\
    + 2\beta_X(\tau)+2\beta_Y(\tau)+ \frac{1}{T^2}\sum_{k,\ell=1}^T\PP{Q_{\alpha'/|\cT|}(\bX,\bY) > Q_{\alpha/|\cT|}(\bX_k,\bY_\ell)}.
\end{multline*}
By definition of $f$, we can also calculate
\begin{align*}
    \frac{1}{T^2}\sum_{k,\ell=1}^T\PP{f(\bX_k,\bY_\ell)\geq Q_{\alpha'/|\cT|}(\bX,\bY)}
    &=\frac{1}{T^2}\sum_{k,\ell=1}^T\PP{\max_{t\in\cT}\psi^*((\bX_k)_t,(\bY_\ell)_t)\geq Q_{\alpha'/|\cT|}(\bX,\bY)}\\
    &\leq \sum_{t\in\cT}\frac{1}{T^2}\sum_{k,\ell=1}^T\PP{\psi^*((\bX_k)_t,(\bY_\ell)_t)\geq Q_{\alpha'/|\cT|}(\bX,\bY)}\\
    &\leq \sum_{t\in\cT}\frac{1}{T^2}\sum_{i,j=1}^T \PP{\psi^*(X_i,Y_j)\geq Q_{\alpha'/|\cT|}(\bX,\bY)}\\
   &\leq \sum_{t\in\cT}\frac{\alpha'}{|\cT|} = \alpha',
\end{align*}
where the last inequality holds by definition of $Q_{\alpha'/|\cT|}$, while the next-to-last inequality holds since, for a fixed choice of $t$, there is a one-to-one correspondence between $(k,\ell)\in[T]\times[T]$ and $(i,j)\in[T]\times[T]$, such that either $(\bX_k)_t=X_i$, $(\bY_\ell)_t=Y_j$ or $\psi^*((\bX_k)_t,(\bY_\ell)_t)=-\infty$.
Next, by~\eqref{eqn:counting},
\[\frac{1}{T^2}\sum_{k,\ell=1}^T \One{t^*\in I_{k,\tau}\cap I_{\ell,\tau}}
= \frac{1}{T^2}\sum_{k,\ell=1}^T \One{(t^*,t^*)\in  I_{k,\tau}\times I_{\ell,\tau}} = \left(\frac{T-\tau}{T}\right)^2 \geq 1-\frac{2\tau}{T},\]
and so $\frac{1}{T^2}\sum_{k,\ell=1}^T \PP{t^*\not\in  I_{k,\tau}\cap I_{\ell,\tau}}\leq \frac{2\tau}{T}$.
Combining our results so far, we have therefore shown
\[
    \PP{\min_{t\in\cT}p_t\leq \frac{\alpha}{|\cT|}}
    \leq \alpha'
    + \frac{2\tau}{T}+ 2\beta_X(\tau)+2\beta_Y(\tau)+ \frac{1}{T^2}\sum_{k,\ell=1}^T\PP{Q_{\alpha'/|\cT|}(\bX,\bY) > Q_{\alpha/|\cT|}(\bX_k,\bY_\ell)}.
\]
Observe that the vector $\bX_k$ is a permutation of the vector $\bX^{T-k-\tau}$ if $k\leq T-\tau-1$, or of $\bX^{2T-k-\tau}$ if $k\geq T-\tau$. 
Similarly, $\bY_\ell$ is a permutation of the vector $\bY^{T-\ell-\tau}$ if $\ell\leq T-\tau-1$, or of $\bY^{2T-\ell-\tau}$ if $\ell\geq T-\tau$. Since the function $Q_a(\bx,\by)$ is invariant to permutations of $\bx$ and $\by$, therefore,
\begin{multline*}\frac{1}{T^2}\sum_{k,\ell=1}^T\PP{Q_{\alpha'/|\cT|}(\bX,\bY) > Q_{\alpha/|\cT|}(\bX_k,\bY_\ell)} 
= \frac{1}{T^2}\sum_{k,\ell=1}^T\PP{Q_{\alpha'/|\cT|}(\bX,\bY) > Q_{\alpha/|\cT|}(\bX^k,\bY^\ell)}\\
= \EE{\frac{1}{T^2}\sum_{k,\ell=1}^T\One{Q_{\alpha'/|\cT|}(\bX,\bY) > Q_{\alpha/|\cT|}(\bX^k,\bY^\ell)}}\leq \begin{cases} 0, & \frac{\alpha'}{|\cT|} \geq \frac{\alpha}{|\cT|} + \frac{2\tau}{T},\\\frac{\alpha'}{\alpha'-\alpha} \cdot \frac{2\tau}{T},&\textnormal{ otherwise},\end{cases}\end{multline*}
where the last step uses the following deterministic bound:
\begin{lemma}\label{lem:counting_argument}
    For any $a<b\in[0,1]$, and any $\bx\in\cX^T$, $\by\in\cY^T$,
    \[\frac{1}{T^2}\sum_{k,\ell=1}^T\One{Q_b(\bx,\by) > Q_a(\bx^k,\by^\ell)} \leq \frac{b}{b-a} \cdot \frac{2\tau}{T}.\]
    Moreover, if $b\geq a + \frac{2\tau}{T}$, then
    \[\frac{1}{T^2}\sum_{k,\ell=1}^T\One{Q_b(\bx,\by) > Q_a(\bx^k,\by^\ell)} =0.\]
\end{lemma}
Combining with our results above, we obtain
\[
    \PP{\min_{t\in\cT}p_t\leq \frac{\alpha}{|\cT|}}
    \leq \alpha'
    + \frac{2\tau}{T}+ 2\beta_X(\tau)+2\beta_Y(\tau)+ \begin{cases} 0, & \frac{\alpha'}{|\cT|} \geq \frac{\alpha}{|\cT|} + \frac{2\tau}{T},\\\frac{\alpha'}{\alpha'-\alpha} \cdot \frac{2\tau}{T},&\textnormal{ otherwise}\end{cases}.
\]
Since this holds for any $\alpha'>\alpha$, by choosing $\alpha' =\alpha + \sqrt{\alpha \cdot \frac{2\tau}{T}}$ we obtain
\[
    \PP{\min_{t\in\cT}p_t\leq \frac{\alpha}{|\cT|}}
    \leq\alpha + \sqrt{\alpha \cdot \frac{8\tau}{T}}
    + \frac{4\tau}{T}+ 2\beta_X(\tau)+2\beta_Y(\tau),
\]
while by choosing $\alpha' = \alpha + \frac{2\tau}{T} \cdot |\cT|$ we obtain
\[
    \PP{\min_{t\in\cT}p_t\leq \frac{\alpha}{|\cT|}}
    \leq\alpha + \frac{2\tau}{T} \cdot |\cT|
    + \frac{2\tau}{T}+ 2\beta_X(\tau)+2\beta_Y(\tau),
\]
which completes the proof.

\subsubsection{Proof for the case $\Delta\geq 1$}
Next we extend the proof to allow for $\Delta\geq 1$. In fact, we will see that this can be reduced
to the case $\Delta=0$.

Define $\tilde{\cX} = \cX^{\Delta+1}$ and $\tilde{\cY} = \cY^{\Delta+1}$, and $\tilde{T} = T-\Delta$. 
Define also $\tilde{X}_i = (X_i,\dots,X_{i+\Delta})$ for $i\in[\tilde{T}]$, and $\tilde\bX = (\tilde{X}_1,\dots,\tilde{X}_{\tilde{T}})$, and similarly
define $\tilde{Y}_j = (Y_j,\dots,Y_{j+\Delta})$ for $j\in[\tilde{T}]$, and $\tilde\bY = (\tilde{Y}_1,\dots,\tilde{Y}_{\tilde{T}})$.
Finally, note that $\psi$ is a function taking inputs in $\tilde\cX\times\tilde\cY$, by definition.

Therefore, the time-shifted p-value $p_t$ defined in~\eqref{eqn:pval_event}, 
using data streams $\bX,\bY$ of dimension $T$ with observation windows of length $\Delta+1$, is exactly
equivalent to running the same procedure on data streams $\tilde\bX,\tilde\bY$ of dimension $\tilde{T}$,
with observation windows of length $1$ (i.e., $\Delta=0$). Applying the conclusion of Theorem~\ref{thm2} for the case
$\Delta=0$ then yields
\[
    \PP{\min_{t\in\cT}p_t\leq \frac{\alpha}{|\cT|}}
    \leq\alpha +  \min \left\{\sqrt{\alpha \cdot \frac{8\tilde\tau}{\tilde{T}}}, \frac{2\tilde{\tau}}{\tilde{T}}\cdot (|\cT -1|) \right\}
    + \frac{4\tilde\tau}{\tilde{T}}+ 2\beta_{\tilde{X}}(\tilde\tau)+2\beta_{\tilde{Y}}(\tilde\tau),
\]
for any $\tilde\tau\in\{0,\dots,\tilde{T}-1\}$,
where we observe that the $\beta$-mixing coefficients are now computed using the distributions
of the data streams $\tilde\bX$ and $\tilde\bY$ rather than $\bX$ and $\bY$.
Now we bound these coefficients.
Let $\bX'$ be an i.i.d.\ copy of $\bX$, and define
$\tilde{X}'_i =  (X'_i,\dots,X'_{i+\Delta})$ for $i\in[\tilde{T}]$; then $\tilde\bX'=(\tilde{X}'_1,\dots,\tilde{X}'_{\tilde{T}})$ is an i.i.d.\ copy of $\tilde\bX$.
 Let $\tilde\tau \geq \Delta$. For any $k$, we need to compare the data subvectors
\[(\tilde{X}_1,\dots,\tilde{X}_k,\tilde{X}_{k+\tilde\tau+1},\dots,\tilde{X}_{\tilde{T}})\textnormal{ and }(\tilde{X}_1,\dots,\tilde{X}_k,\tilde{X}'_{k+\tilde\tau+1},\dots,\tilde{X}'_{\tilde{T}}).\]
By definition, these are equal to
\[\big( (X_1,\dots,X_{\Delta+1}),\dots,(X_k,\dots,X_{\Delta+k}),(X_{k+\tilde\tau+1},\dots,X_{\Delta+k+\tilde\tau+1}),\dots,(X_{T-\Delta},\dots,X_T)\big)\]
and
\[\big( (X_1,\dots,X_{\Delta+1}),\dots,(X_k,\dots,X_{\Delta+k}),(X'_{k+\tilde\tau+1},\dots,X'_{\Delta+k+\tilde\tau+1}),\dots,(X'_{T-\Delta},\dots,X'_T)\big).\]
Therefore,
\begin{multline*}\dtv\big((\tilde{X}_1,\dots,\tilde{X}_k,\tilde{X}_{k+\tilde\tau+1},\dots,\tilde{X}_{\tilde{T}}),(\tilde{X}_1,\dots,\tilde{X}_k,\tilde{X}'_{k+\tilde\tau+1},\dots,\tilde{X}'_{\tilde{T}})\big)\\
\leq \dtv\big( (X_1,\dots,X_{\Delta+k},X_{k+\tilde\tau+1},\dots,X_T) , (X_1,\dots,X_{\Delta+k},X'_{k+\tilde\tau+1},\dots,X'_T)
\leq \beta_X(\tilde\tau-\Delta),\end{multline*}
where the last step holds by definition of the $\beta$-mixing coefficients for $\bX$.
Therefore, we have shown that
\[\beta_{\tilde{X}}(\tilde\tau) \leq \beta_X(\tilde\tau-\Delta)\]
for all $\tilde\tau\in\{\Delta,\dots,T-\Delta-1\}$. An identical calculation shows that 
\[\beta_{\tilde{Y}}(\tilde\tau) \leq \beta_Y(\tilde\tau-\Delta).\]
Therefore, plugging this in to our calculations above
yields
\[
    \PP{\min_{t\in\cT}p_t\leq \frac{\alpha}{|\cT|}}
    \leq\alpha +  \min \left\{\sqrt{\alpha \cdot \frac{8\tilde\tau}{\tilde{T}}}, \frac{2\tilde{\tau}}{\tilde{T}}\cdot (|\cT -1|) \right\}
    + \frac{4\tilde\tau}{\tilde{T}}+ 2\beta_X(\tilde\tau-\Delta)+2\beta_Y(\tilde\tau-\Delta),
\]
for any  $\tilde\tau\in\{\Delta,\dots,T-\Delta-1\}$. For any $\tau\leq T-2\Delta-1$, this completes the proof when we choose $\tilde\tau = \tau+\Delta$ and plug in $\tilde{T}=T-\Delta$.
If instead $\tau\geq T-2\Delta$ then the result of the theorem holds trivially.

\subsection{Proof of Theorem~\ref{thm3}}

Define a function $P_c$, mapping data streams to a value in $[0,1]$, as
\[P_c(\bx,\by) = 
\frac{1}{T^2} \sum_{i,j=1}^T  \One{ \psi(\bx_{[i]},\by_{[j]}) + c \geq\psi(\bx,\by)},\]
for any $c\geq 0$, where
we define 
the ``wrap-around'' data stream
\[\bx_{[i]} = (x_i,\dots,x_T,x_1,\dots,x_{i-1}),\]
and similarly for $\by_{[j]}$. Note that
\[p^{+\varepsilon} = P_{\varepsilon}(\bX,\bY),\]
by definition. Moreover, by construction, for any $\bx,\by$ and any $c\geq 0$ we must have
\begin{equation}\label{eqn:thm3_deterministic_sum}\frac{1}{T^2}\sum_{i,j=1}^T\One{P_c(\bx_{[i]},\by_{[j]}) \leq \alpha}\leq \alpha.\end{equation}
In addition, below we will prove that, for all $i,j\in[T]$ and any $c\geq 0$,
\begin{equation}\label{eqn:thm3_rotate_ineq}
\PP{P_{c+4\gamma(\tau)}(\bX,\bY) \leq \alpha}\leq \PP{P_c(\bX_{[i]},\bY_{[j]}) \leq \alpha} + 2\beta_X(\tau)+2\beta_Y(\tau).\end{equation}
Combining these two bounds (with $c=\varepsilon-4\gamma(\tau)\geq 0$), we then have
\begin{align*}
    \PP{p^{+\varepsilon}\leq \alpha}
    &=\PP{P_{\varepsilon}(\bX,\bY)\leq \alpha}\\
    &\leq \frac{1}{T^2}\sum_{i,j=1}^T\PP{P_{\varepsilon-4\gamma(\tau)}(\bX_{[i]},\bY_{[j]})\leq \alpha}+ 2\beta_X(\tau)+2\beta_Y(\tau)\textnormal{ by~\eqref{eqn:thm3_rotate_ineq}}\\
    &=\EE{\frac{1}{T^2}\sum_{i,j=1}^T\One{P_{\varepsilon-4\gamma(\tau)}(\bX_{[i]},\bY_{[j]})\leq \alpha}}+ 2\beta_X(\tau)+2\beta_Y(\tau)\\
    &\leq \alpha +  2\beta_X(\tau)+2\beta_Y(\tau)\textnormal{ by~\eqref{eqn:thm3_deterministic_sum}}.
\end{align*}

To complete the proof, we now need to verify~\eqref{eqn:thm3_rotate_ineq}.
We will need the following lemma, which applies the stability assumption to establish properties of the function $P_c$.
\begin{lemma}\label{lem_for_thm3}
Fix any $k,k',\ell,\ell'\in[T]$, and define
\[P^+_c(\bx,\by) = \sup_{\bx',\by'}\left\{P_c(\bx',\by') : \bx_{I_{k,\tau}\cap I_{k',\tau}} = \bx'_{I_{k,\tau}\cap I_{k',\tau}}, \ \by_{I_{\ell,\tau}\cap I_{\ell',\tau}} = \by'_{I_{\ell,\tau}\cap I_{\ell',\tau}}\right\}.\]
If $\psi$ satisfies Assumption~\ref{asm:stability}, then
for any $c\geq 0$ and any $\bx,\by$,
it holds that
\[P_c(\bx,\by) \leq P^+_c(\bx,\by) \leq P_{c+4\gamma(\tau)}(\bx,\by).\]
\end{lemma}
Now fix any $i,j\in[T]$, and define $k=T-i+1$, $\ell=T-j+1$, $k'=\ell'=T$. Then we have
\begin{align*}
    &\PP{P_{c+4\gamma(\tau)}(\bX,\bY) \leq \alpha}\\
    &\leq \PP{P^+_{c}(\bX,\bY) \leq \alpha}\\
    &\leq \PP{P^+_{c}(\bX_{[i]},\bY_{[j]}) \leq \alpha} \\
    &\hspace{.5in}{}+ \dtv\Big(\big(\bX_{I_{T-i+1,\tau}\cap I_{T,\tau}},\bY_{I_{T-j+1,\tau}\cap I_{T,\tau}}\big), \big((\bX_{[i]})_{I_{T-i+1,\tau}\cap I_{T,\tau}},(\bY_{[j]})_{I_{T-j+1,\tau}\cap I_{T,\tau}}\big)\Big)\\
    &\leq \PP{P_c(\bX_{[i]},\bY_{[j]}) \leq \alpha} \\
    &\hspace{.5in}{}+ \dtv\Big(\big(\bX_{I_{T-i+1,\tau}\cap I_{T,\tau}},\bY_{I_{T-j+1,\tau}\cap I_{T,\tau}}\big), \big((\bX_{[i]})_{I_{T-i+1,\tau}\cap I_{T,\tau}},(\bY_{[j]})_{I_{T-j+1,\tau}\cap I_{T,\tau}}\big)\Big).
\end{align*}
Here the first and last step hold by Lemma~\ref{lem_for_thm3} (since $\psi$ satisfies Assumption~\ref{asm:stability}), while the middle step holds because $P^+_c(\bx,\by)$ depends on $(\bx,\by)$ only through $\bx_{I_{T-i+1,\tau}\cap I_{T,\tau}}$ and $\by_{I_{T-j+1,\tau}\cap I_{T,\tau}}$, for any $c$. Finally, applying~\eqref{eqn:dtv_for_thm3} to bound the total variation distance, we have completed the proof of~\eqref{eqn:thm3_rotate_ineq}.

\subsection{Proofs of lemmas}
\subsubsection{Proof of Lemma~\ref{lem:dtv_rotate}}
The statements for $\bX$ and $\bY$ are identical, so we only prove the bound for $\bX$. We split into cases:
\paragraph{Case 1: $1 \leq i\leq \tau+1$.}
Then
\[I_{T-i+1,\tau}\cap I_{T,\tau} = \{\tau+1,\dots,T-i+1\}.\]
Therefore, 
\[\bX_{I_{T-i+1,\tau}\cap I_{T,\tau}} 
= (X_{\tau+1},\dots,X_{T-i+1}),\]
and
\[(\bX_{[i]})_{I_{T-i+1,\tau}\cap I_{T,\tau}} 
= (X_{\tau+i},\dots,X_T).\]
Since $\bX$ is stationary, these two vectors are equal in distribution.

\paragraph{Case 2: $\tau + 2 \leq i\leq T-\tau$.}
We have
\[I_{T-i+1,\tau}\cap I_{T,\tau} = \{\tau+1,\dots,T-i+1,T-i+\tau+2,\dots,T\}.\]
Therefore, 
\[\bX_{I_{T-i+1,\tau}\cap I_{T,\tau}} 
= (X_{\tau+1},\dots,X_{T-i+1},X_{T-i+\tau+2},\dots,X_T),\]
and
\[(\bX_{[i]})_{I_{T-i+1,\tau}\cap I_{T,\tau}} 
= (X_{\tau+i},\dots,X_T,X_{\tau+1},\dots,X_{i-1}).\]
Now let $\bX'$ be an i.i.d.\ copy of $\bX$. Then, by definition of $\beta$-mixing,
\[\dtv\big((X_{\tau+1},\dots,X_{T-i+1},X_{T-i+\tau+2},\dots,X_T),
(X_{\tau+1},\dots,X_{T-i+1},X_{T-i+\tau+2}',\dots,X_T')\big)\leq \beta_X(\tau),\]
and
\[\dtv\big((X_{\tau+i},\dots,X_T,X_{\tau+1},\dots,X_{i-1}),(X_{\tau+i},\dots,X_T,X_{\tau+1}',\dots,X_{i-1}')\big)\leq \beta_X(\tau).\]
Since $\bX,\bX'$ are stationary and i.i.d., it also holds that 
\[(X_{\tau+1},\dots,X_{T-i+1},X_{T-i+\tau+2}',\dots,X_T') \eqd (X_{\tau+i},\dots,X_T,X_{\tau+1}',\dots,X_{i-1}'). \]
By the triangle inequality, then, we have proved the desired bound.

\paragraph{Case 3: $T-\tau+1\leq i\leq T$ and $i \ge \tau +2$.} Then
\[I_{T-i+1,\tau}\cap I_{T,\tau} = \{T-i+\tau+2,\dots,T\}.\]
Therefore, 
\[\bX_{I_{T-i+1,\tau}\cap I_{T,\tau}} 
= (X_{T-i+\tau+2},\dots,X_T),\]
and
\[(\bX_{[i]})_{I_{T-i+1,\tau}\cap I_{T,\tau}} 
= (X_{\tau+1},\dots,X_{i-1}).\]
Since $\bX$ is stationary, these two vectors are equal in distribution.

\subsubsection{Proof of Lemma~\ref{lem:counting_argument}}
Fix any $k,\ell$ and suppose $Q_b(\bx,\by) > Q_a(\bx^k,\by^\ell)$. 
By definition of the function $Q_a$, we must therefore have
\[\frac{1}{T^2}\sum_{i,j=1}^T \One{\psi^*((\bx^k)_i,(\by^\ell)_j) \geq Q_b(\bx,\by)} < a.\]
On the other hand, again by definition,
\[\frac{1}{T^2}\sum_{i,j=1}^T \One{\psi^*(x_i,y_j) \geq Q_b(\bx,\by)} \geq b.\]
By construction of $\bx^k,\by^\ell$,
\[\psi^*((\bx^k)_i,(\by^\ell)_j) = \begin{cases} \psi^*(x_i,y_j) , &(i,j)\in I_{k,\tau}\times I_{\ell,\tau},\\
-\infty, & \textnormal{ otherwise}.\end{cases}\]
Therefore,
\[\sum_{i,j=1}^T \One{\psi^*((\bx^k)_i,(\by^\ell)_j) \geq Q_b(\bx,\by)}
= \sum_{i,j=1}^T \One{\psi^*(x_i,y_j) \geq Q_b(\bx,\by)} \cdot\One{(i,j)\in I_{k,\tau}\times I_{\ell,\tau}}.\]
Combining these calculations, then, we must have
\[\frac{1}{T^2}\sum_{i,j=1}^T \One{\psi^*(x_i,y_j) \geq Q_b(\bx,\by)} \cdot\One{(i,j)\not\in I_{k,\tau}\times I_{\ell,\tau}}
>  b-a.
\]
If $b \ge a+ \frac{2\tau}{T}$, then 
\begin{multline*}
    \frac{2\tau}{T} \le b-a < \frac{1}{T^2}\sum_{i,j=1}^T \One{\psi^*(x_i,y_j) \geq Q_b(\bx,\by)} \cdot\One{(i,j)\not\in I_{k,\tau}\times I_{\ell,\tau}} \\ \le \frac{1}{T^2}\sum_{i,j=1}^T\One{(i,j)\not\in I_{k,\tau}\times I_{\ell,\tau}} \le \frac{2\tau}{T},
\end{multline*}
which is a contradiction. Thus,
\[\One{\frac{1}{T^2}\sum_{i,j=1}^T \One{\psi^*(x_i,y_j) \geq Q_b(\bx,\by)} \cdot\One{(i,j)\not\in I_{k,\tau}\times I_{\ell,\tau}}
>  b-a} = 0,\]
when $b-a \ge \frac{2\tau}{T}$. Therefore,
\begin{multline*}
    \frac{1}{T^2}\sum_{k,\ell=1}^T\One{Q_b(\bx,\by) > Q_a(\bx_k,\by_\ell)} \\ \le \frac{1}{T^2}\sum_{k,\ell=1}^T\One{\frac{1}{T^2}\sum_{i,j=1}^T \One{\psi^*(x_i,y_j) \geq Q_b(\bx,\by)} \cdot\One{(i,j)\not\in I_{k,\tau}\times I_{\ell,\tau}}>  b-a} = 0.
\end{multline*}

For the general case, where we do not assume $b\geq a + \frac{2\tau}{T}$, we instead calculate that
\begin{align*}
    &\frac{1}{T^2}\sum_{k,\ell=1}^T\One{Q_b(\bx,\by) > Q_a(\bx_k,\by_\ell)}\\
    &\leq \frac{1}{T^2}\sum_{k,\ell=1}^T\One{\frac{1}{T^2}\sum_{i,j=1}^T \One{\psi^*(x_i,y_j) \geq Q_b(\bx,\by)} \cdot\One{(i,j)\not\in I_{k,\tau}\times I_{\ell,\tau}}
>  b-a}\\
    &\leq \frac{1}{b-a} \cdot \frac{1}{T^2}\sum_{k,\ell=1}^T \left(\frac{1}{T^2}\sum_{i,j=1}^T \One{\psi^*(x_i,y_j) \geq Q_b(\bx,\by)} \cdot\One{(i,j)\not\in I_{k,\tau}\times I_{\ell,\tau}}\right)\\
    &= \frac{1}{b-a} \cdot \frac{1}{T^2}\sum_{i,j=1}^T  \One{\psi^*(x_i,y_j) \geq Q_b(\bx,\by)}
    \left(\frac{1}{T^2}\sum_{k,\ell=1}^T \One{(i,j)\not\in I_{k,\tau}\times I_{\ell,\tau}}\right)\\
    &\leq \frac{1}{b-a} \cdot \frac{1}{T^2}\sum_{i,j=1}^T  \One{\psi^*(x_i,y_j) \geq Q_b(\bx,\by)}\cdot \frac{2\tau}{T}\\
    &\leq \frac{1}{b-a} \cdot b\cdot \frac{2\tau}{T},
\end{align*}
which completes the proof.
\subsubsection{Proof of Lemma~\ref{lem_for_thm3}}
    The lower bound on $P^+_c(\bx,\by)$ is trivial, since we can simply take $\bx'=\bx$ and $\by'=\by$. Now we prove the upper bound. Fix any $\bx',\by'$ with $\bx_{I_{k,\tau}\cap I_{k',\tau}} = \bx'_{I_{k,\tau}\cap I_{k',\tau}}$ and $\by_{I_{\ell,\tau}\cap I_{\ell',\tau}} = \by'_{I_{\ell,\tau}\cap I_{\ell',\tau}}$.
    Let $\bx''$ be defined as:
    \[x''_i = \begin{cases} x_i, &i\in I_{k,\tau}, \\ x'_i, &i\not\in I_{k,\tau}.\end{cases}\]
    Then by construction, $\bx''_{I_{k,\tau}}=\bx_{I_{k,\tau}}$ and $\bx''_{I_{k',\tau}} = \bx'_{I_{k',\tau}}$. Consequently, $\bx''\stackrel{\tau}{\asymp}\bx$ and $\bx''\stackrel{\tau}{\asymp}\bx'$. We can similarly construct a vector $\by''$ such that $\by''\stackrel{\tau}{\asymp}\by$ and $\by''\stackrel{\tau}{\asymp}\by'$. Then we have
    \[\left|\psi(\bx,\by) - \psi(\bx',\by')\right|
    \leq \left|\psi(\bx,\by) - \psi(\bx'',\by'')\right| + \left|\psi(\bx'',\by'') - \psi(\bx',\by')\right| \leq 2\gamma(\tau),\]
    where for the last step we apply Assumption~\ref{asm:stability} to each term. Moreover, by definition of the relation $\stackrel{\tau}{\asymp}$, for any $\mathbf{z},\mathbf{z}'$ and any $k\in[T]$, 
    \[\mathbf{z}\stackrel{\tau}{\asymp}\mathbf{z}' \ \Longleftrightarrow \ \mathbf{z}_{[k]}\stackrel{\tau}{\asymp}\mathbf{z}'_{[k]}.\]
    Consequently,
    for any $i,j\in[T]$, it holds that $\bx''_{[i]}\stackrel{\tau}{\asymp}\bx_{[i]}$ and $\bx''_{[i]}\stackrel{\tau}{\asymp}\bx'_{[i]}$, and similarly, $\by''_{[j]}\stackrel{\tau}{\asymp}\by_{[j]}$ and $\by''_{[j]}\stackrel{\tau}{\asymp}\by_{[j]}'$.
    Therefore,
    \[\left|\psi(\bx_{[i]},\by_{[j]}) - \psi(\bx'_{[i]},\by'_{[j]})\right|
    \leq \left|\psi(\bx_{[i]},\by_{[j]}) - \psi(\bx''_{[i]},\by''_{[j]})\right| + \left|\psi(\bx''_{[i]},\by''_{[j]}) - \psi(\bx'_{[i]},\by'_{[j]})\right| \leq 2\gamma(\tau),\]
    where again for the last step we apply Assumption~\ref{asm:stability} to each term. Therefore, for any $i,j$,
    \[\textnormal{If $\psi(\bx'_{[i]},\by'_{[j]}) + c \geq \psi(\bx',\by')$ then 
    $\psi(\bx_{[i]},\by_{[j]}) + c + 4\gamma(\tau) \geq \psi(\bx,\by)$}.\]
    Therefore,
    \begin{multline*}P_c(\bx',\by') = \frac{1}{T^2}\sum_{i,j=1}^T \One{\psi(\bx'_{[i]},\by'_{[j]}) + c \geq \psi(\bx',\by')}\\\leq \frac{1}{T^2}\sum_{i,j=1}^T \One{\psi(\bx_{[i]},\by_{[j]}) + c + 4\gamma(\tau) \geq \psi(\bx,\by)} = P_{c+4\gamma(\tau)}(\bx,\by).\end{multline*}
Since this is true for all $\bx',\by'$ satisfying 
$\bx_{I_{k,\tau}\cap I_{k',\tau}} = \bx'_{I_{k,\tau}\cap I_{k',\tau}}$ and $\by_{I_{\ell,\tau}\cap I_{\ell',\tau}} = \by'_{I_{\ell,\tau}\cap I_{\ell',\tau}}$, we therefore have
$P^+_c(\bx,\by)\leq P_{c+4\gamma(\tau)}(\bx,\by)$, as desired.

\section{An example for matching the inflation factor in Theorem~\ref{thm1}} \label{app:thm1_factor}
In the simulation experiments in Section~\ref{sec:experiments}, the p-values for testing the midpoint in time (i.e., $p_t$, for testing whether an event occurred at time $t \approx T/2$) result in a false positive probability that is approximately bounded by $\alpha$, empirically.
In contrast, Theorem~\ref{thm1} only establishes an upper bound of $\approx 4\alpha$ on the probability of a false positive. This raises the question of whether the theoretical bound is loose. However, we will show with examples that the inflation factor $1/\mathsf{mar}(t)$ in Theorem~\ref{thm1} is tight (when $\alpha$ is sufficiently small). In particular, when $t \approx T/2$, the factor of $4$ cannot be improved. For simplicity, we will work in the setting $\Delta=0$, and will not consider the multiple testing problem (that is, we only test for an event at a single time $t$).

\begin{proposition}
    Fix any $T\geq 1$, any $\alpha \in (0,1)$, and any $t \in \{1,\cdots,T\}$. Assume $\sqrt{\alpha}T$ is an integer. Then there exist independent and stationary time series $\bX$, $\bY \in \R^T$, and a function $\psi: \R \times \R \to \R$, such that
    \[\PP{p_t \le \alpha} \ge \frac{\alpha}{\mathsf{mar}(t)} \cdot  \left(\frac{1}{1+\sqrt{\frac{\alpha}{\mathsf{mar}(t)}}}\right)^2.\]
\end{proposition}

Note that when $\alpha \ll \mathsf{mar}(t)$, this lower bound is approximately equal to $\frac{\alpha}{\mathsf{mar}(t)}$, which matches the upper bound shown in Theorem~\ref{thm1}.

\begin{proof}
    Without loss of generality, we may assume $t\leq \frac{T+1}{2}$, so that $\mathsf{mar}(t) = (T/t)^2$. Let $m = \sqrt{\alpha} T$ and let $N$ be a large positive integer with $N>T$ (the value of $N$ will be specified below).

    Define a time series $\bX$ by taking the product of two components.
    \begin{itemize}
        \item A cyclic component: let $I_X \sim \mathrm{Uniform} \left(\{0, \cdots, m+t-2\}  \right)$. Define
        \[\bX^{\textnormal{cyclic}} = (X_1^{\textnormal{cyclic}}, \cdots, X_T^{\textnormal{cyclic}}), \ \textnormal{ where } X_i^{\textnormal{cyclic}} = \One{\mathrm{mod} (I_X+i, m+t-1) \notin \{1, \cdots, t-1 \}}.\]
        This component is equivalent to sampling $T$ consecutive entries from the infinite periodic string
        \[ \cdots, \underbrace{1, \cdots, 1}_{m \textnormal{ times}}, \underbrace{0,\cdots, 0}_{t-1 \textnormal{ times}},  \underbrace{1, \cdots, 1}_{m \textnormal{ times}}, \underbrace{0,\cdots, 0}_{t-1 \textnormal{ times}}, \cdots \]

        \item A component that is monotonic with high probability: let $J_X \sim \mathrm{Uniform} \left(\{0, \cdots, N-1\}  \right)$. Define
        \[\bX^{\textnormal{mono}} = (X_1^{\textnormal{mono}}, \cdots, X_T^{\textnormal{mono}}), \ \textnormal{ where } X_i^{\textnormal{mono}} = N - \mathrm{mod}(J_X+i, N).\]
        This component is equivalent to sampling $T$ consecutive entries from the infinite periodic string
        \[ \cdots, N, N-1,  \cdots,2, 1, N, N-1, \cdots, 2, 1, \cdots \]
    \end{itemize}
    Let $\bX = (X_1^{\textnormal{cyclic}} \cdot X_1^{\textnormal{mono}}, \cdots, X_T^{\textnormal{cyclic}} \cdot X_T^{\textnormal{mono}})$, the elementwise product of the two components $\bX^{\textnormal{cyclic}}$ and $\bX^{\textnormal{mono}}$, with each of these two components generated independently. Since each component has a stationary distribution by construction, $\bX$ is a stationary time series.

    Next, let $\bY$ be generated independently from the same distribution as $\bX$: we independently sample $I_Y\sim \mathrm{Uniform} \left(\{0, \cdots, m+t-2\}  \right)$ and $J_Y\sim\mathrm{Uniform} \left(\{0, \cdots, N-1\}  \right)$ independently, and define $\bY^{\textnormal{cyclic}}$, $\bY^{\textnormal{mono}}$, and $\bY$ analogously. 
    
    Let $\psi: \R \times \R \to \R$ be defined as $\psi(x,y) = (x+y)\cdot\one{x> 0}\cdot\one{y>0}$. Consider the event that
    \begin{equation}\label{eqn:event_for_lower_bound}I_X, I_Y \in \{0, \cdots,m-1 \} \textnormal{ and } J_X, J_Y \in \{0,\cdots, N-T-1\}.\end{equation}
    On this event, it holds that
    \[\mathrm{mod} (I_X+t, m+t-1) \in \{t,t+1,\dots,m+t-2,0\},\]
    which is the complement of $\{1,\cdots,t-1\}$.
    Thus, $X^{\textnormal{cyclic}}_t = 1$. Analogously, $Y^{\textnormal{cyclic}}_t = 1$.
    Meanwhile, since $N>T$, it holds that
    \[ X_i^{\textnormal{mono}} = N - J_X -i, \textnormal{ and }  Y_j^{\textnormal{mono}} = N - J_Y - j\]
    for all $i,j \in [T]$.

    Therefore, on the event in~\eqref{eqn:event_for_lower_bound}, $\Psi_t =\psi(X_t, Y_t) = X_t + Y_t = 2N-J_X-J_Y - 2t >0$, and for each $i,j \in [T]$ we can calculate
    \begin{itemize}
        \item If $X^{\textnormal{cyclic}}_i = 0$ or $Y^{\textnormal{cyclic}}_j = 0$, then $\Psi_{i,j} = 0$, and so
        \[\One{\Psi_{i,j} \ge \Psi_t} = \One{0 \ge \Psi_t} = 0; \]

        \item If $X^{\textnormal{cyclic}}_i=Y^{\textnormal{cyclic}}_j=1$, then $\Psi_{i,j} = X^{\textnormal{mono}}_i + Y^{\textnormal{mono}}_j$, and so
        \[\One{\Psi_{i,j} \ge \Psi_t} = \One{N-J_X-i + N-J_Y - j \ge 2N-J_X-J_Y - 2t} = \One{i+j \le 2t}.\]
    \end{itemize}
    Therefore, on the event in~\eqref{eqn:event_for_lower_bound},
    \begin{align*}
        p_t &=\frac{1}{T^2}\sum_{i,j=1}^T \One{\Psi_{i,j}\geq\Psi_t}
        = \frac{1}{T^2} \sum_{i,j=1}^T \One{X^{\textnormal{cyclic}}_i= Y^{\textnormal{cyclic}}_j=1, i+j \le 2t} \\
        &\leq \frac{1}{T^2} \sum_{i,j=1}^{T} \One{X^{\textnormal{cyclic}}_i= Y^{\textnormal{cyclic}}_j=1, i\le 2t-1,j \le 2t-1}\\
        &=\frac{1}{T^2}\left(\sum_{i=1}^{2t-1}\One{X^{\textnormal{cyclic}}_i=1}\right)\left(\sum_{j=1}^{2t-1}\One{Y^{\textnormal{cyclic}}_j=1}\right).\end{align*}
    Since $X^{\textnormal{cyclic}}_t = 1$ by~\eqref{eqn:event_for_lower_bound}, by definition of $\bX^{\textnormal{cyclic}}$ we can verify that $\sum_{i=1}^{2t-1}\One{X^{\textnormal{cyclic}}_i=1} \leq m$, and the analogous bound holds for $\bY$ (see Figure~\ref{fig:illustrate_circ1} for an illustration of this calculation). Therefore, on the event defined in~\eqref{eqn:event_for_lower_bound},
    \[p_t \leq \frac{m^2}{T^2}=\alpha.\]

    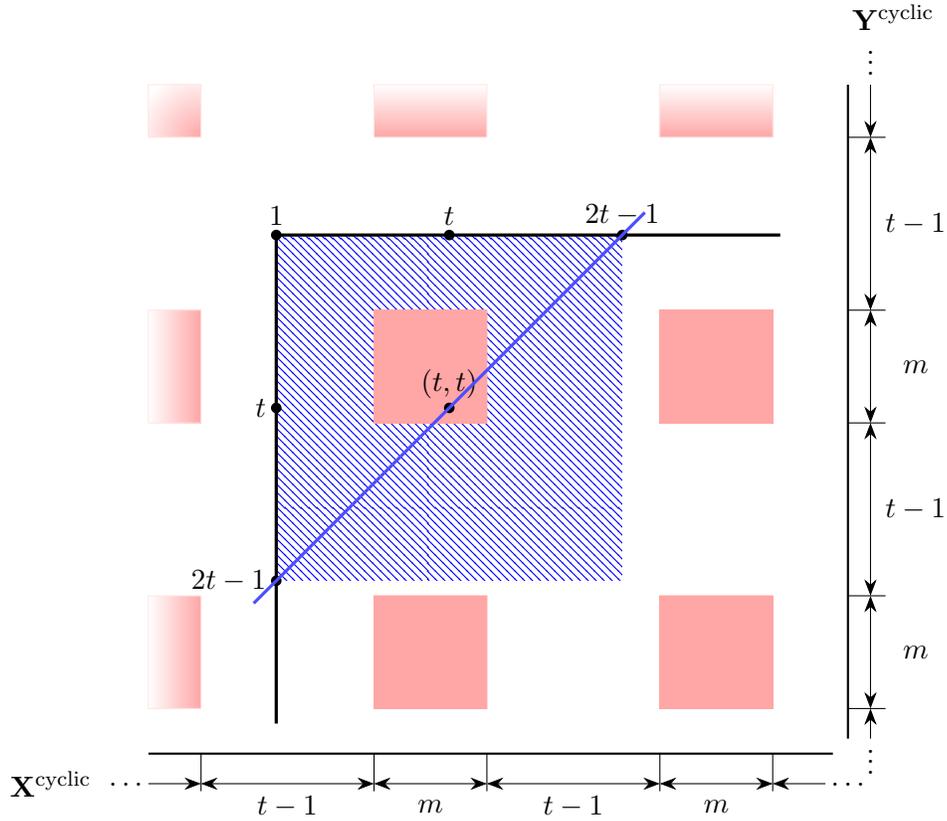
\begin{figure}
  \centering
  \begin{tikzpicture}[
  x=1cm,y=1cm,
    >={Stealth[length=2.2mm]}
    ]

\def\W{8}
\def\yC{5.2}

\fill[pattern=north west lines, pattern color=blue] (0.7,\yC+0.2) rectangle (5.3,\yC-4.4);

\shade[
  left color=red!35,
  right color=red!35,
  draw=red!35,
  line width=0.1pt
] (2, \yC-2.3) rectangle (3.5, \yC-0.8);

\shade[
  left color=red!35,
  right color=red!35,
  draw=red!35,
  line width=0.1pt
] (5.8,\yC-2.3) rectangle (7.3, \yC-0.8);
\shade[
  left color=red!35,
  right color=red!35,
  draw=red!35,
  line width=0.1pt
] (2,\yC-6.1) rectangle (3.5,\yC-4.6);
\shade[
  left color=red!35,
  right color=red!35,
  draw=red!35,
  line width=0.1pt
] (5.8,\yC-6.1) rectangle (7.3,\yC-4.6);

\shade[
  left color=red!0,
  right color=red!35,
  draw=red!10,
  line width=0.1pt
] (-1, \yC-2.3) rectangle (-0.3, \yC-0.8);

\shade[
  left color=red!0,
  right color=red!35,
  draw=red!10,
  line width=0.1pt
] (-1,\yC-6.1) rectangle (-0.3,\yC-4.6);

\shade[
  top color=red!0,
  bottom color=red!35,
  draw=red!10,
  line width=0.1pt
] (2, \yC+1.5) rectangle (3.5,\yC+2.2);

\shade[
  top color=red!0,
  bottom color=red!35,
  draw=red!10,
  line width=0.1pt
] (5.8,\yC+1.5) rectangle (7.3,\yC+2.2);

\begin{scope}
  \clip (-1,\yC+1.5) rectangle (-0.3,\yC+2.2);
  \shade[inner color=red!35, outer color=red!0]
        (-0.3,\yC+1.5) circle [radius=0.9];  
\end{scope}
\draw[red!10] (-1,\yC+1.5) rectangle (-0.3,\yC+2.2);

\draw[line width=1.2pt, black]
  (0.7,\yC+0.2) -- (\W-0.6, \yC+0.2);
\draw[line width=1.2pt, black]
  (0.7, \yC+0.2) -- (0.7, \yC+1.7-\W);

\fill[black] (0.7,\yC+0.2) circle (2pt);
\fill[black] (3, \yC+0.2) circle (2pt);
\fill[black] (0.7, \yC-2.1) circle (2pt);
\fill[black] (5.3, \yC+0.2) circle (2pt);
\fill[black] (0.7, \yC-4.4) circle (2pt);
\node[above] at (3, \yC+0.2) {$t$};
\node[above] at (5.3, \yC+0.2) {$2t-1$};
\node[above] at (0.7,\yC+0.2) {$1$};
\node[left] at (0.7, \yC-2.1) {$t$};
\node[left] at (0.7, \yC-4.4) {$2t-1$};
\fill[black] (3, \yC-2.1) circle (2pt);
\node[above] at (3, \yC-2.1) {$ (t,t)$};

\draw[blue!70, line width=1.2pt]
  (0.4, \yC-4.7) -- (5.6, \yC+0.5);

\draw[line width=0.8pt, black]
  (-1,-1.5) -- (\W+0.1, -1.5);
\draw[line width=0.8pt, black]
  (\W+0.3, -1.3) -- (\W+0.3, \yC+2.2);

\foreach \x in {-0.3, 2,3.5, 5.8, 7.3}{
  \draw (\x,-1.5) -- ++(0,-0.5);
}
\foreach \y in {\yC+1.5, \yC-2.3, \yC-0.8, \yC-6.1, \yC-4.6}{
  \draw (\W+0.3,\y) -- ++(0.5,0);
}

\node at (-2.3, -1.9) {$\bX^{\textnormal{cyclic}}$};
\draw[->, thin] (-1,-1.9) -- (-0.3,-1.9);
\draw[<->, thin] (-0.3,-1.9) -- (2,-1.9);
\node at (0.85, -2.2) {$t-1$};
\draw[<->, thin] (2,-1.9) -- (3.5,-1.9);
\node at (2.75, -2.2) {$m$};
\draw[<->, thin] (3.5,-1.9) -- (5.8,-1.9);
\node at (4.65, -2.2) {$t-1$};
\draw[<->, thin] (5.8,-1.9) -- (7.3,-1.9);
\node at (6.55, -2.2) {$m$};
\draw[<-, thin] (7.3,-1.9) -- (\W,-1.9);
\node at (\W+0.3,-1.9) {$\hdots$};
\node at (-1.3,-1.9) {$\hdots$};

\node at (\W+0.9, \yC+3.1) {$\bY^{\textnormal{cyclic}}$};
\draw[->, thin] (\W+0.6, \yC+2.2) -- (\W+0.6,\yC+1.5);
\draw[<->, thin] (\W+0.6,\yC+1.5) -- (\W+0.6, \yC-0.8);
\node at (\W+1.2, \yC+0.35) {$t-1$};
\draw[<->, thin] (\W+0.6,\yC-2.3) -- (\W+0.6, \yC-0.8);
\node at (\W+1.2, \yC-1.55) {$m$};
\draw[<->, thin] (\W+0.6,\yC-2.3) -- (\W+0.6, \yC-4.6);
\node at (\W+1.2, \yC-3.45) {$t-1$};
\draw[<->, thin] (\W+0.6,\yC-4.6) -- (\W+0.6, \yC-6.1);
\node at (\W+1.2, \yC-5.35) {$m$};
\draw[<-, thin] (\W+0.6, \yC-6.1) -- (\W+0.6,\yC-6.5);
\node at (\W+0.6,\yC-6.7) {$\vdots$};
\node at (\W+0.6,\yC+2.6) {$\vdots$};

\end{tikzpicture}
\caption{An illustration of $\One{X^{\textnormal{cyclic}}_i = Y^{\textnormal{cyclic}}_j =1} $, over all integers $i,j$, with red shading indicating a value of $1$. Note that on the event in~\eqref{eqn:event_for_lower_bound}, it holds that $X^{\textnormal{cyclic}}_t=Y^{\textnormal{cyclic}}_t=1$ (i.e., the point $(t,t)$ lies in a $m\times m$ red shaded square, as in the illustration). On this event, no other red shaded squares overlap with the region $(i,j)\in\{1,\dots,2t-1\}\times \{1,\dots,2t-1\}$ (indicated by shading with blue lines).}
\label{fig:illustrate_circ1}
\end{figure}
   
    Finally, by definition of the distributions of $I_X,J_X,I_Y,J_Y$, the event in~\eqref{eqn:event_for_lower_bound} holds with probability
    \begin{multline*}\left(\frac{m}{m+t-1} \right)^2 \cdot \left(\frac{N-T}{N} \right)^2 = \frac{\alpha T^2}{\left(t + \sqrt{\alpha}T - 1 \right)^2} \cdot \left(\frac{N-T}{N} \right)^2\\ = \frac{\alpha}{\left(\frac{t}{T} + \sqrt{\alpha} - \frac{1}{T} \right)^2} \cdot \left(1 - \frac{T}{N}\right)^2
    = \frac{\alpha}{\mathsf{mar}(t)}\cdot\left(\frac{1 - \frac{T}{N} }{1 + \sqrt{\frac{\alpha}{\mathsf{mar}(t)}}- \frac{1}{t}}\right)^2,\end{multline*}
    and therefore,
    \[\PP{p_t\leq \alpha}\geq \frac{\alpha}{\mathsf{mar}(t)}\cdot\left(\frac{1 - \frac{T}{N} }{1 + \sqrt{\frac{\alpha}{\mathsf{mar}(t)}}- \frac{1}{t}}\right)^2.\]
    Finally, by choosing a sufficiently large positive integer $N$, we can ensure that
    \[\frac{1 - \frac{T}{N} }{1 + \sqrt{\frac{\alpha}{\mathsf{mar}(t)}}- \frac{1}{t}} \geq \frac{1}{1 + \sqrt{\frac{\alpha}{\mathsf{mar}(t)}}},\] which completes the proof.
\end{proof}

\bibliographystyle{plainnat}
\bibliography{bib}

\end{document}